\documentclass[AMA,STIX1COL]{WileyNJD-v2}

\usepackage{bbm}
\usepackage{graphicx}
\usepackage{wrapfig}
\usepackage{caption}
\usepackage{subcaption}
\usepackage[section]{placeins} 

\usepackage{fancyhdr} 
\usepackage{lastpage} 
\usepackage{extramarks} 

\usepackage{xcolor} 
\usepackage{enumerate} 
\usepackage{paralist} 
\usepackage{amsmath}
\usepackage{amsthm}
\usepackage{amssymb, mathtools}
\usepackage{pifont, stmaryrd} 
\usepackage{etoolbox} 
\usepackage[T1]{fontenc}
\usepackage{bibentry} 
\makeatletter\let\saved@bibitem\@bibitem\makeatother 
\makeatletter\let\@bibitem\saved@bibitem\makeatother 

\usepackage{bm} 
\usepackage{amsfonts} 

\DeclareMathOperator*{\argmin}{arg\,min}

\newcommand{\func}[3]{\ensuremath{#1 : #2 \rightarrow #3}}

\newcommand{\norm}[1]{\ensuremath{\left\| #1 \right\|}}

\newcommand{\pder}[2]{\ensuremath{\frac{\partial #1}{\partial #2}}}

\newcommand{\Dcal}{\ensuremath{\mathcal{D}}}

\newcommand{\Tcal}{\ensuremath{\mathcal{T}}}

\newcommand{\Rbb}{\ensuremath{\mathbb{R} }}

\usepackage{tikz}
\usepackage{pgfplots}
\usepackage{pgfplotstable, filecontents, booktabs}
\pgfplotsset{compat=1.9}

\usetikzlibrary{pgfplots.groupplots}
\usepgfplotslibrary{fillbetween}
\usetikzlibrary{calc,fit,matrix,arrows,automata,positioning,shapes}
\usetikzlibrary{arrows.meta}

\pgfplotsset{select coords between index/.style 2 args={
    x filter/.code={
        \ifnum\coordindex<#1\fi
        \ifnum\coordindex>#2\fi
    }
}}

\tikzset{
 invisible/.style={opacity=0},
 visible on/.style={alt={#1{}{invisible}}},
 alt/.code args={<#1>#2#3}{%
   \alt<#1>{\pgfkeysalso{#2}}{\pgfkeysalso{#3}}
 },
}

\newcommand{\colorbarMatlabParula}[5]{
\begin{tikzpicture}
\begin{axis}[
   hide axis, scale only axis,
   height=0pt, width=0pt,
   colormap={parula}{rgb255=(62,38,168) rgb255=(62,39,172) rgb255=(63,40,175) rgb255=(63,41,178) rgb255=(64,42,180) rgb255=(64,43,183) rgb255=(65,44,186) rgb255=(65,45,189) rgb255=(66,46,191) rgb255=(66,47,194) rgb255=(67,48,197) rgb255=(67,49,200) rgb255=(67,50,202) rgb255=(68,51,205) rgb255=(68,52,208) rgb255=(69,53,210) rgb255=(69,55,213) rgb255=(69,56,215) rgb255=(70,57,217) rgb255=(70,58,220) rgb255=(70,59,222) rgb255=(70,61,224) rgb255=(71,62,225) rgb255=(71,63,227) rgb255=(71,65,229) rgb255=(71,66,230) rgb255=(71,68,232) rgb255=(71,69,233) rgb255=(71,70,235) rgb255=(72,72,236) rgb255=(72,73,237) rgb255=(72,75,238) rgb255=(72,76,240) rgb255=(72,78,241) rgb255=(72,79,242) rgb255=(72,80,243) rgb255=(72,82,244) rgb255=(72,83,245) rgb255=(72,84,246) rgb255=(71,86,247) rgb255=(71,87,247) rgb255=(71,89,248) rgb255=(71,90,249) rgb255=(71,91,250) rgb255=(71,93,250) rgb255=(70,94,251) rgb255=(70,96,251) rgb255=(70,97,252) rgb255=(69,98,252) rgb255=(69,100,253) rgb255=(68,101,253) rgb255=(67,103,253) rgb255=(67,104,254) rgb255=(66,106,254) rgb255=(65,107,254) rgb255=(64,109,254) rgb255=(63,110,255) rgb255=(62,112,255) rgb255=(60,113,255) rgb255=(59,115,255) rgb255=(57,116,255) rgb255=(56,118,254) rgb255=(54,119,254) rgb255=(53,121,253) rgb255=(51,122,253) rgb255=(50,124,252) rgb255=(49,125,252) rgb255=(48,127,251) rgb255=(47,128,250) rgb255=(47,130,250) rgb255=(46,131,249) rgb255=(46,132,248) rgb255=(46,134,248) rgb255=(46,135,247) rgb255=(45,136,246) rgb255=(45,138,245) rgb255=(45,139,244) rgb255=(45,140,243) rgb255=(45,142,242) rgb255=(44,143,241) rgb255=(44,144,240) rgb255=(43,145,239) rgb255=(42,147,238) rgb255=(41,148,237) rgb255=(40,149,236) rgb255=(39,151,235) rgb255=(39,152,234) rgb255=(38,153,233) rgb255=(38,154,232) rgb255=(37,155,232) rgb255=(37,156,231) rgb255=(36,158,230) rgb255=(36,159,229) rgb255=(35,160,229) rgb255=(35,161,228) rgb255=(34,162,228) rgb255=(33,163,227) rgb255=(32,165,227) rgb255=(31,166,226) rgb255=(30,167,225) rgb255=(29,168,225) rgb255=(29,169,224) rgb255=(28,170,223) rgb255=(27,171,222) rgb255=(26,172,221) rgb255=(25,173,220) rgb255=(23,174,218) rgb255=(22,175,217) rgb255=(20,176,216) rgb255=(18,177,214) rgb255=(16,178,213) rgb255=(14,179,212) rgb255=(11,179,210) rgb255=(8,180,209) rgb255=(6,181,207) rgb255=(4,182,206) rgb255=(2,183,204) rgb255=(1,183,202) rgb255=(0,184,201) rgb255=(0,185,199) rgb255=(0,186,198) rgb255=(1,186,196) rgb255=(2,187,194) rgb255=(4,187,193) rgb255=(6,188,191) rgb255=(9,189,189) rgb255=(13,189,188) rgb255=(16,190,186) rgb255=(20,190,184) rgb255=(23,191,182) rgb255=(26,192,181) rgb255=(29,192,179) rgb255=(32,193,177) rgb255=(35,193,175) rgb255=(37,194,174) rgb255=(39,194,172) rgb255=(41,195,170) rgb255=(43,195,168) rgb255=(44,196,166) rgb255=(46,196,165) rgb255=(47,197,163) rgb255=(49,197,161) rgb255=(50,198,159) rgb255=(51,199,157) rgb255=(53,199,155) rgb255=(54,200,153) rgb255=(56,200,150) rgb255=(57,201,148) rgb255=(59,201,146) rgb255=(61,202,144) rgb255=(64,202,141) rgb255=(66,202,139) rgb255=(69,203,137) rgb255=(72,203,134) rgb255=(75,203,132) rgb255=(78,204,129) rgb255=(81,204,127) rgb255=(84,204,124) rgb255=(87,204,122) rgb255=(90,204,119) rgb255=(94,205,116) rgb255=(97,205,114) rgb255=(100,205,111) rgb255=(103,205,108) rgb255=(107,205,105) rgb255=(110,205,102) rgb255=(114,205,100) rgb255=(118,204,97) rgb255=(121,204,94) rgb255=(125,204,91) rgb255=(129,204,89) rgb255=(132,204,86) rgb255=(136,203,83) rgb255=(139,203,81) rgb255=(143,203,78) rgb255=(147,202,75) rgb255=(150,202,72) rgb255=(154,201,70) rgb255=(157,201,67) rgb255=(161,200,64) rgb255=(164,200,62) rgb255=(167,199,59) rgb255=(171,199,57) rgb255=(174,198,55) rgb255=(178,198,53) rgb255=(181,197,51) rgb255=(184,196,49) rgb255=(187,196,47) rgb255=(190,195,45) rgb255=(194,195,44) rgb255=(197,194,42) rgb255=(200,193,41) rgb255=(203,193,40) rgb255=(206,192,39) rgb255=(208,191,39) rgb255=(211,191,39) rgb255=(214,190,39) rgb255=(217,190,40) rgb255=(219,189,40) rgb255=(222,188,41) rgb255=(225,188,42) rgb255=(227,188,43) rgb255=(230,187,45) rgb255=(232,187,46) rgb255=(234,186,48) rgb255=(236,186,50) rgb255=(239,186,53) rgb255=(241,186,55) rgb255=(243,186,57) rgb255=(245,186,59) rgb255=(247,186,61) rgb255=(249,186,62) rgb255=(251,187,62) rgb255=(252,188,62) rgb255=(254,189,61) rgb255=(254,190,60) rgb255=(254,192,59) rgb255=(254,193,58) rgb255=(254,194,57) rgb255=(254,196,56) rgb255=(254,197,55) rgb255=(254,199,53) rgb255=(254,200,52) rgb255=(254,202,51) rgb255=(253,203,50) rgb255=(253,205,49) rgb255=(253,206,49) rgb255=(252,208,48) rgb255=(251,210,47) rgb255=(251,211,46) rgb255=(250,213,46) rgb255=(249,214,45) rgb255=(249,216,44) rgb255=(248,217,43) rgb255=(247,219,42) rgb255=(247,221,42) rgb255=(246,222,41) rgb255=(246,224,40) rgb255=(245,225,40) rgb255=(245,227,39) rgb255=(245,229,38) rgb255=(245,230,38) rgb255=(245,232,37) rgb255=(245,233,36) rgb255=(245,235,35) rgb255=(245,236,34) rgb255=(245,238,33) rgb255=(246,239,32) rgb255=(246,241,31) rgb255=(246,242,30) rgb255=(247,244,28) rgb255=(247,245,27) rgb255=(248,247,26) rgb255=(248,248,24) rgb255=(249,249,22) rgb255=(249,251,21) },
   colorbar horizontal,
   point meta min=#1, point meta max=#5,
   colorbar style={width=10cm, xtick={#1,#2,#3,#4,#5}}]

\addplot [draw=none] coordinates {(0,0)};
\end{axis}
\end{tikzpicture}
}

\articletype{Article Type}%

\received{TBD}
\revised{TBD}
\accepted{TBD}

\begin{document}




\title{Accurate quantification of blood flow wall shear stress using simulation-based imaging: a synthetic, comparative study}

\author[1]{Charles J. Naudet}

\author[2]{Johannes T\"{o}ger}

\author[1]{Matthew J. Zahr*}

\authormark{C. J. Naudet \textsc{et al}}

\address[1]{\orgdiv{Department of Aerospace and Mechanical Engineering}, \orgname{University of Notre Dame}, \orgaddress{\state{Indiana}, \country{USA}}}

\address[2]{\orgdiv{Department of Clinical Sciences Lund, Clinical Physiology}, \orgname{Lund University, Sk\r{a}ne University Hospital}, \orgaddress{\state{Lund}, \country{Sweden}}}

\corres{*Matthew J. Zahr, 300B Cushing Hall, Notre Dame, IN 46656. \email{mzahr@nd.edu}}


\abstract[Abstract]{Simulation-based imaging (SBI) is a blood flow imaging technique that optimally
fits a computational fluid dynamics (CFD) simulation to low-resolution,
noisy magnetic resonance (MR) flow data to produce a high-resolution
velocity field. In this work, we study the effectivity of SBI in
predicting wall shear stress (WSS) relative to standard magnetic
resonance imaging (MRI) postprocessing techniques using two
synthetic numerical experiments: flow through an idealized,
two-dimensional stenotic vessel and a model of an adult aorta.
In particular,
we study the sensitivity of these two approaches with
respect to the Reynolds number of the underlying flow,
the resolution of the MRI data, and the noise in the MRI
data. We found that the SBI WSS reconstruction:
1) is insensitive to Reynolds number over the range considered
   ($\mathrm{Re} \leq 1000$),
2) improves as the amount of MRI data increases and provides
   accurate reconstructions with as little as three MRI voxels per
   diameter, and
3) degrades linearly as the noise in the data increases with a
   slope determined by the resolution of the MRI data.
We also consider the sensitivity of SBI to the
resolution of the CFD mesh and found there is flexibility
in the mesh used for SBI, although the WSS reconstruction
becomes more sensitive to other parameters, particularly
the resolution of the MRI data, for coarser meshes. This
indicates a fundamental trade-off between scan time
(i.e., MRI data quality and resolution) and reconstruction time
using SBI, which is inherently different than the trade-off
between scan time and reconstruction quality observed
in standard MRI postprocessing techniques.}


\jnlcitation{\cname{%
\author{C.J. Naudet}, 
\author{J. T\"{o}ger}, 
\author{M.J. Zahr}} (\cyear{2021}), 
\ctitle{Accurate quantification of blood flow wall shear stress using simulation-based imaging: a synthetic, comparative study}, \cjournal{International Journal of Numerical Methods in Biomedical Engineering}, \cvol{TBD}.}

\maketitle

\section{Background}

Magnetic resonance imaging (MRI) is a powerful method to investigate
cardiovascular physiology. High-resolution \textit{in vivo} images can
help understand patient-specific blood flow and provide important quantitative
biomarkers such as wall shear stress (WSS). However, these
methods are limited by a fundamental trade-off between scan time, resolution,
and noise \cite{edelstein1986intrinsic,portnoy2009information}, which limits their utility for
applications that demand very high-resolution images such as infants and
children with congenital heart disease \cite{markl_advanced_2016}.
This is further complicated when estimates of biomarkers must be extracted
from these flow images with poorly resolved features. This motivates the
need for imaging methods that can use sparse, noisy data to provide sufficiently
high-resolution reconstructions that can be used to accurately compute
quantitative biomarkers.


Recent advances in compressed sensing
\cite{compressed_sensing_holland,4Dflow_undersamp,compressed_sensing_lustig,parallel_imaging_tariq}
and machine learning \cite{lundervold2019overview,vishnevskiy2020deep} have
been used to enhance MRI-based flow reconstructions, leading to improved
image quality and reduced scan times. Neural networks have proven effective
in taking sparse or missing data in k-space and accurately reconstructing
them into natural images. Machine learning approaches are valuable, but face
drawbacks of high training cost and not being patient-specific.
Several simulation-based imaging (SBI) methods
\cite{de2016temporal,de20144d,funke2019variational,gaidzik2019transient,goenezen20164d,morbiducci2010outflow,rispoli2015using,toger_zahr_4Dflow,wood2001combined}
exist that match a computational fluid dynamics (CFD) simulation to magnetic
resonance (MR) flow data by optimizing free parameters of the CFD simulation (usually
boundary and initial conditions), use a metric to measure the difference
between CFD and MR flow imaging, and update the free parameters to minimize
the cost function via their own strategy. The method used in this work
\cite{toger_zahr_4Dflow} uses a high-order CFD discretization, efficient
adjoint-based PDE-constrained optimization, and a novel objective function
that mimics the point-spread function of MRI scanners. This method was shown to
effectively reconstruct very high-resolution velocity fields from limited
MRI data and match ground truth values for both a controlled water tank
experiment and an \textit{in vivo} clinical application. However, the
quality to which the method predicts quantitative biomarkers has not
been considered to date.

In this paper, we study the accuracy to which the SBI approach of
\cite{toger_zahr_4Dflow} predicts the WSS distribution relative
to standard MRI postprocessing approaches using two synthetic
numerical experiments. We focus on WSS because it is known to correlate
to atherosclerosis, the formation and rupture aneurysms, as well as
numerous congenital heart diseases
\cite{wss_imp_boussel,wss_motiv_cheng,wss_motiv_Fedak,wss_imp_frydrychowicz,wss_imp_groen,wss_motiv_Lasheras,atheros_wss_malek,wss_motiv_papaioannaou,wss_imp_reneman,wss_motiv_rosenthal,wss_imp_shaaban}. Furthermore,
it has proven difficult to estimate directly using standard
MRI postprocessing techniques
\cite{van2013wall,potters2014measuring,potters2015volumetric}
because WSS requires accurate estimation of the velocity gradient
that can be difficult using only piecewise constant voxel data.
Current methods for postprocessing MR flow data to obtain quantities
of interest, e.g., phase-contrast (PC) MRI velocity mapping, Fourier
velocity encoding (FVE), and intravoxel velocity standard deviation
mapping \cite{quant_wss,quant_wss_numerical,quant_wss_fem}, can be
unreliable.

We study the impact of the Reynolds number
of the underlying flow and the resolution and noise of the MRI data
to understand the sensitivity of each method with respect to these
critical parameters. Noise is a critical source of error for \textit{in vivo}
imaging and extraction of biomarkers \cite{ha2016hemodynamic} and becomes
increasingly problematic as the resolution of the MRI grid increases or if
faster scans are required, e.g. for sedated children or to reduce costs of
health care. On the other hand, the resolution of the MRI data can lead to
higher resolution images and more accurate biomarker computations; however,
it is usually accompanied with increased noise and requires longer scans.
The Reynolds number is studied because it has been observed that the accuracy
of the WSS computed directly from MRI data decreases as the Reynolds number
of the flow increases \cite{quant_wss_numerical}. Furthermore, we study the
impact of the resolution of the CFD mesh used in SBI because this determines
the overall cost of the reconstruction.

\section{Methods}
\label{sec:methods}
Two synthetic numerical experiments were conducted to compare the accuracy of MRI methods with SBI in measuring wall shear stress. To perform comprehensive studies, SBI was simplified for use in two-dimensional, time-independent problems,
instead of three-dimensional, unsteady problems as in
our previous study \cite{toger_zahr_4Dflow}.
The experiments are synthetic in the sense that no \textit{in vivo} MRI flow
data or geometries were used; synthetic data was constructed to be
representative of a realistic situation and consistent with \textit{in vivo}
measurements to the extent possible, e.g., regarding noise levels, MRI
resolution, and extraction of MRI data from a velocity field. The synthetic
experiments allow for a highly controlled study with a known reference
(``truth'') flow so the impact of various parameters, e.g., Reynolds number,
noise, MRI grid resolution, on WSS reconstruction
accuracy can be isolated and identified. The remainder
of this section will describe the numerical experiments in detail.
Section~\ref{sec:methods:synthetic} introduces the setup of the synthetic
experiments, Section~\ref{sec:methods:mri} details the creation of synthetic
MRI data and WSS reconstruction, and Section~\ref{sec:methods:sbi} reviews the
SBI method and corresponding WSS computation.

\subsection{Synthetic experiments}
\label{sec:methods:synthetic}
\begin{figure}
 \centering
 \raisebox{-0.5\height}{\input{_py/stenosis0_geom.tikz}}
 \raisebox{-0.5\height}{\begin{tikzpicture}
\begin{axis}[
axis equal image,
axis x line*=bottom,
axis y line*=left,
width=0.65\textwidth,
xtick={0, 2.5, 6.6, 10.4, 12.5},
ytick={-7.9, 0, 4.9, 11.6},
ymax=12.73,
xmax=13.35,
xticklabel style={/pgf/number format/precision=2, /pgf/number format/fixed},
xmin=-0.87,
ymin=-8.78,
yticklabel style={/pgf/number format/precision=2, /pgf/number format/fixed}]
\addplot []
graphics [xmin=-0.31,xmax=12.79,ymin=-8.22,ymax=12.17] { 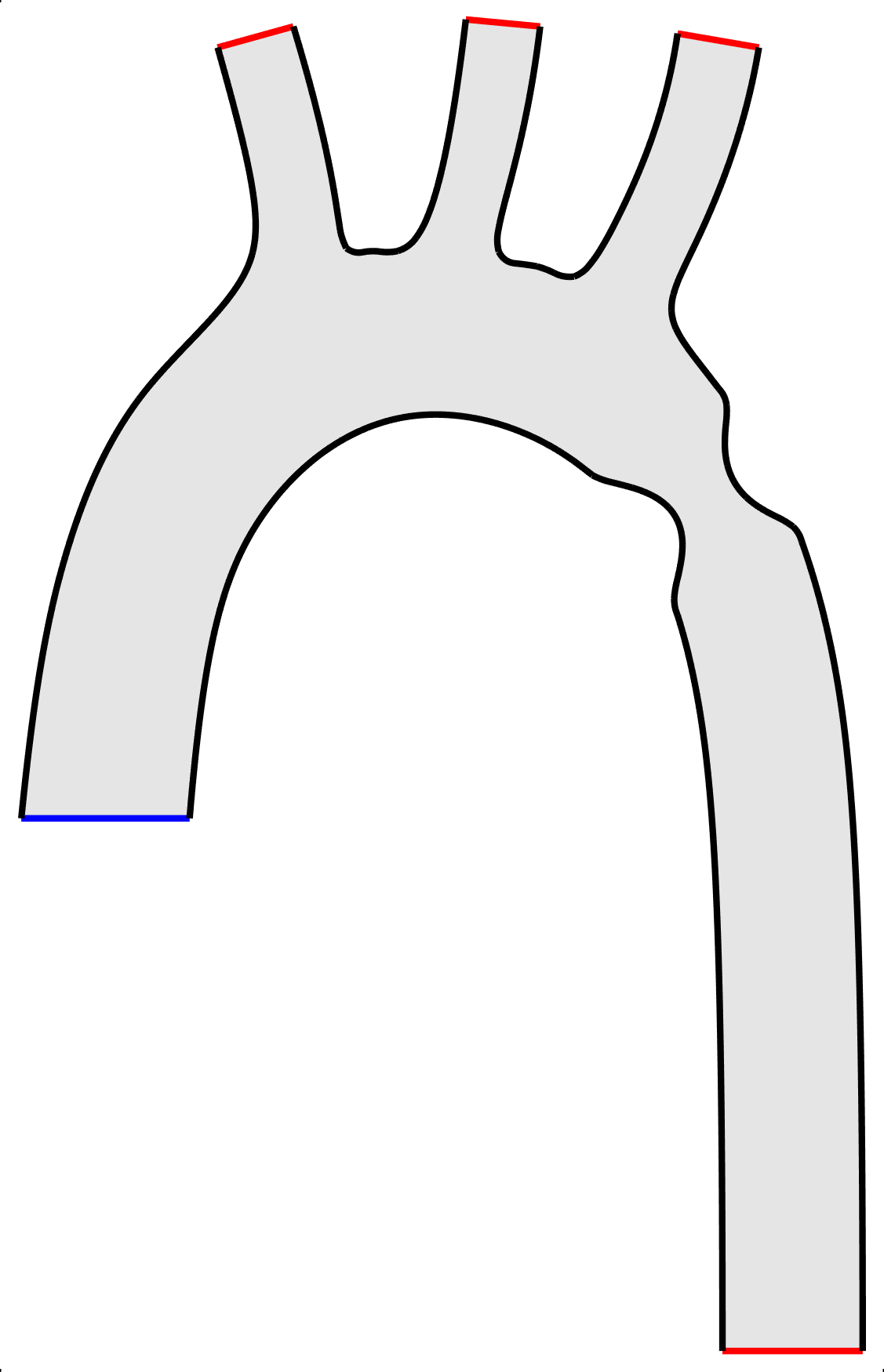};

\addplot [magenta, dashed, thick, forget plot]
coordinates {
( 0.00000000e+00, -7.90000000e+00)
( 1.30000000e+01, -7.90000000e+00)
( 1.30000000e+01,  1.16000000e+01)
( 0.00000000e+00,  1.16000000e+01)
( 0.00000000e+00, -7.90000000e+00)};

\end{axis}
\end{tikzpicture}}
 \caption{Idealized stenotic vessel with 60\% grade (\textit{left}) and
  coarctated aorta with 40\% grade (\textit{right}) geometry $\Omega$ used
  for numerical experiments with boundaries: $\partial\Omega_\text{w}$
  (\ref{line:stenosis0-wall}),
  $\partial\Omega_\text{in}$ (\ref{line:stenosis0-inflow}),
  $\partial\Omega_\text{out}$ (\ref{line:stenosis0-outflow}).
  The MRI scan region for each geometry is indicated with
  (\ref{line:stenosis0-mribox}). Units: cm.}
 \label{fig:geom}
\end{figure}

We use two synthetic experiments to investigate the accuracy of WSS
measurements from SBI relative to standard MRI methods: 1) flow through an
idealized stenotic vessel and 2) flow through an idealized aorta with a
coarctation.
The geometry of the stenotic vessel with 60\% grade is the set
$\Omega\subset\Rbb^d$ ($d=2$ in this work; formulas not pertaining
to domain geometries hold for general $d$), defined as
\begin{equation} \label{eqn:geom}
    \Omega \coloneqq \{(x, \pm y(x)) \mid x \in [0, 6]\}, \quad y : x \mapsto B_0 - \frac{A}{\sqrt{2\pi\sigma^2}}\text{exp}\left(\frac{(x-c)^2}{2\sigma^2}\right)
\end{equation}
where $B_0=0.3$ cm, $c=3$ cm, $\sigma=0.6$ cm, $A=0.18$ cm$^2$
(Figure~\ref{fig:geom}).
The geometry of the coarctated aorta (40\% grade) is shown in
Figure~\ref{fig:geom}.
In practice, these geometries would be obtained from MRI scans; however, we
chose to explicitly define the geometry to ensure a controlled setting in
which to study WSS reconstruction. Most of the studies in this work are
conducted using the vessel due to its simplicity; the aorta is used to
confirm our findings on a more realistic geometry.

The blood flow is modeled as an incompressible, Newtonian fluid governed by the Navier-Stokes equations
 \begin{equation} \label{eqn:ins}
     \pder{v}{t} + (v\cdot\nabla)v - \nu\nabla^2 v + \frac{1}{\rho_0}\nabla P = 0, \quad \nabla \cdot v = 0 \quad\text{in}~~\Omega,
 \end{equation}
 where $\rho_0\in\Rbb_{>0}$ is the density of the fluid, $\nu\in\Rbb_{>0}$
 is the kinematic viscosity of the fluid, and $\func{v}{\Omega}{\Rbb^d}$ and
 $\func{P}{\Omega}{\Rbb_{>0}}$ are the
 velocity and pressure, respectively,
 of the fluid implicitly defined as
 the solution of (\ref{eqn:ins}).
 In this work, we assume the fluid is blood and take the material properties
 to be $\rho_0 = 1060$ kg/m$^3$ \cite{trudnowski1974specific} and $\nu = 2.83 \times 10^{-6}$ m$^2$/s \cite{nader2019blood} for
 both test cases.
 We consider the case where the flow has reached a steady state and the time derivative, $\pder{v}{t}$, vanishes. 
 Boundary conditions for the boundaries identified in Figure~\ref{fig:geom} are
\begin{equation}
    v = 0~~\text{on}~~\partial\Omega_\text{w}, \quad
    v = v_\text{in}~~\text{on}~~\partial\Omega_\text{in}, \quad
    \sigma \cdot n = 0~~\text{on}~~\partial\Omega_\text{out}
\end{equation}
where $v_\text{in} : x \mapsto \bar{v}_\mathrm{in} (B_0^2-x_2^2)/B_0^2$ is the
parabolic inlet profile with peak value $\bar{v}_\mathrm{in}\in\Rbb^d$ (for the
stenotic vessel) and $\func{n}{\partial\Omega}{\Rbb^d}$ is the outward unit
normal to the boundary of the domain. The rate-of-strain,
$\func{\epsilon}{\Omega}{\Rbb^{d\times d}}$, and stress,
$\func{\sigma}{\Omega}{\Rbb^{d\times d}}$, tensors are defined as
\begin{equation}
    \epsilon = \frac{1}{2}(\nabla v+\nabla v^T), \qquad
    \sigma = 2 \mu \epsilon + P I_d,
\end{equation}
where $\mu=\rho_0\nu$ is the dynamic viscosity and $I_d$ is the $d\times d$ identity matrix.
The Reynolds number, $\mathrm{Re}\in\Rbb_{>0}$,
 of the flow is defined based on
 the full cross-sectional diameter
 of the geometry, $D\in\Rbb_{>0}$,
 and the peak inlet velocity as
 $\mathrm{Re} = \frac{D \norm{\bar{v}_\mathrm{in}}}{\nu}$.
 The WSS, $\func{\sigma^\mathrm{wss}}{\partial\Omega}{\Rbb}$,
 the quantitative biomarker considered in this work, is defined
 as the magnitude of the tangential component of the surface traction
 \cite{arzani2016characterizations}
\begin{equation} \label{eqn:wss}
    t = \sigma\cdot n, \qquad
    \tau = t - (t\cdot n)n, \qquad
    \sigma^\mathrm{wss} = \norm{\tau},
\end{equation}
where $\func{t}{\partial\Omega}{\Rbb^d}$ is the surface traction and $\func{\tau}{\partial\Omega}{\Rbb^d}$ is its
tangential component.

\begin{figure}
 \centering
 \begin{tikzpicture}
\begin{groupplot} [
group style={group size = 2 by 3, horizontal sep = 0.4cm, vertical sep = 0.8cm},
title style={at={(current bounding box.north west)}, anchor=west}]
\nextgroupplot[axis equal image, axis lines=none, width=0.58\textwidth, ymax=0.09, xmax=1, title=(a) Reference mesh, xmin=0, ymin=-0.09]
\addplot []
graphics [xmin=0,xmax=1,ymin=-0.05,ymax=0.05] { 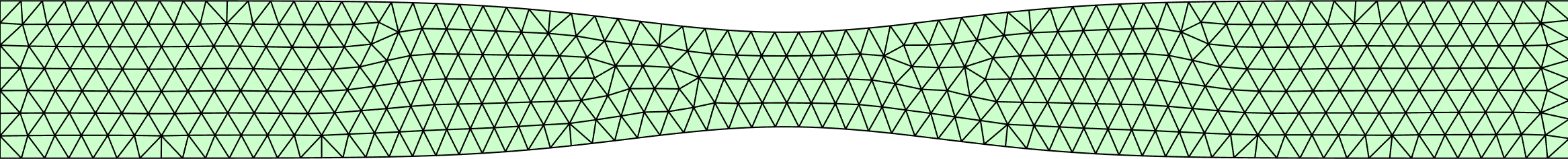};

\nextgroupplot[axis equal image, axis lines=none, width=0.58\textwidth, ymax=0.09, xmax=1, title=(b) Reference velocity, xmin=0, ymin=-0.09]
\addplot []
graphics [xmin=0,xmax=1,ymin=-0.05,ymax=0.05] { 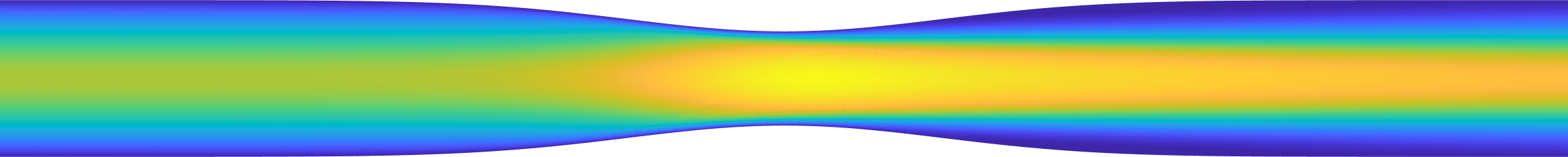};

\nextgroupplot[axis equal image, axis lines=none, width=0.58\textwidth, ymax=0.09, xmax=1, title=(c) MR data (no noise), xmin=0, ymin=-0.09]
\addplot []
graphics [xmin=0,xmax=1,ymin=-0.09,ymax=0.09] { 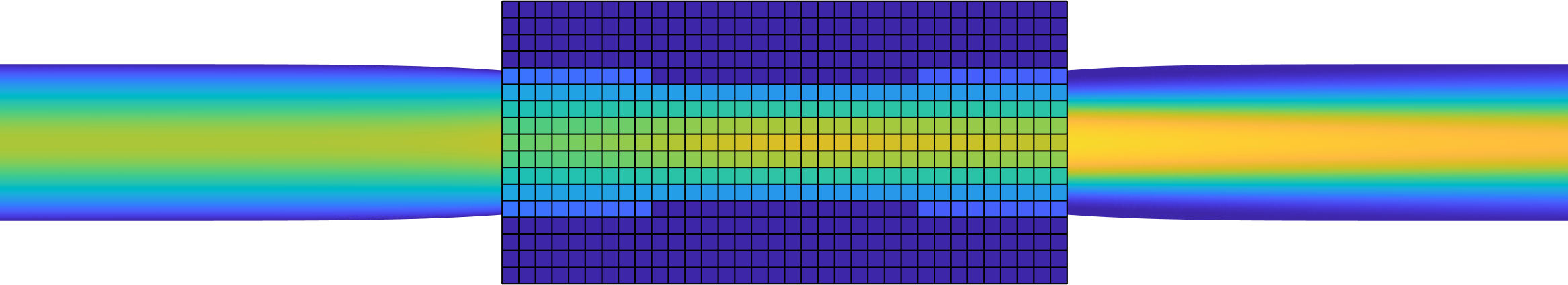};

\nextgroupplot[axis equal image, axis lines=none, width=0.58\textwidth, ymax=0.09, xmax=1, title=(d) MRI data (noise), xmin=0, ymin=-0.09]
\addplot []
graphics [xmin=0,xmax=1,ymin=-0.09,ymax=0.09] { 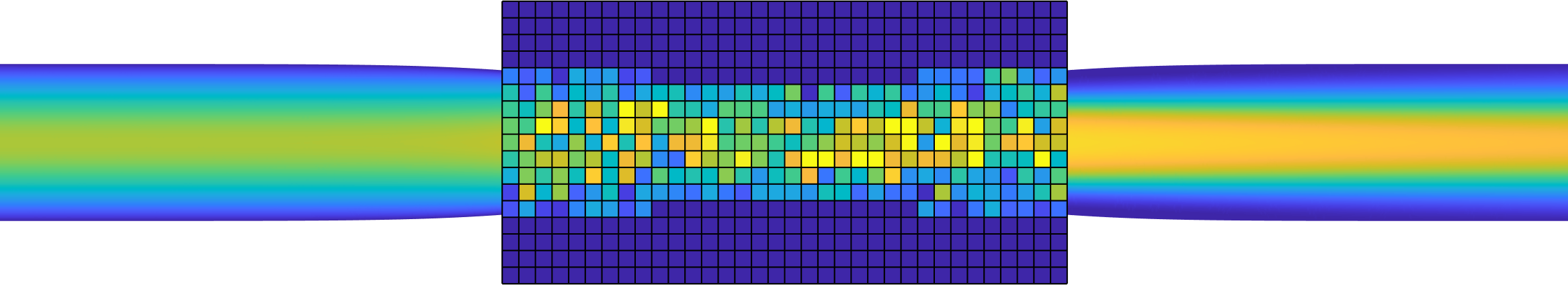};

\nextgroupplot[axis equal image, axis lines=none, width=0.58\textwidth, ymax=0.09, xmax=1, title=(e) SBI reconstructed velocity, xmin=0, ymin=-0.09]
\addplot []
graphics [xmin=0,xmax=1,ymin=-0.05,ymax=0.05] { 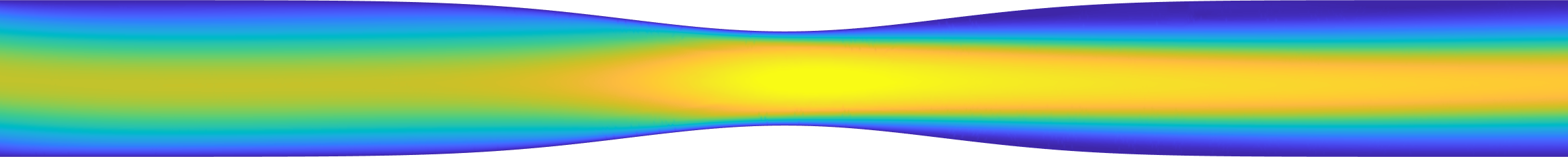};

\nextgroupplot[axis equal image, axis lines=none, width=0.58\textwidth, ymax=0.09, xmax=1, title=(f) MR data from SBI, xmin=0, ymin=-0.09]
\addplot []
graphics [xmin=0,xmax=1,ymin=-0.09,ymax=0.09] { 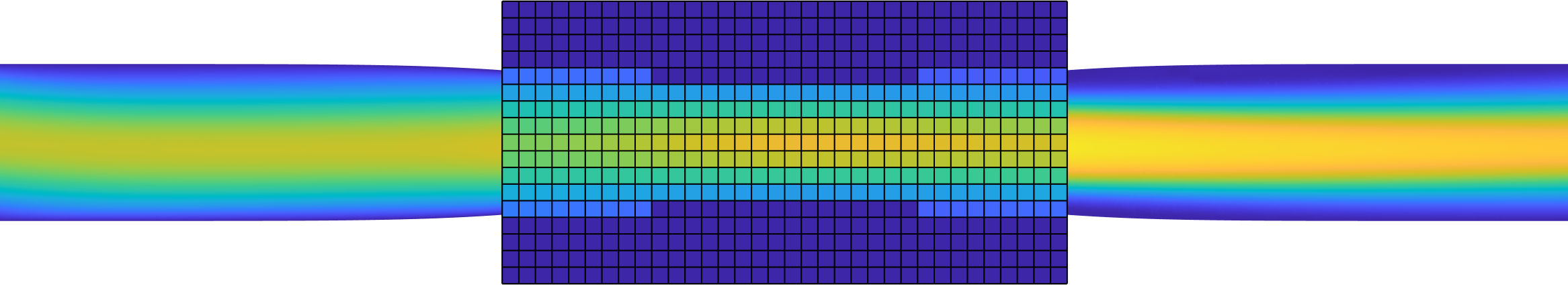};

\end{groupplot}\end{tikzpicture}
 \colorbarMatlabParula{0}{3.5}{7}{10.5}{14}
 \caption{Computational mesh (a) and corresponding velocity
 field (b) used to define the true flow ($v_h$) through the stenotic
 vessel ($\mathrm{Re}=1000$). The velocity field is mapped to the MR data
 space ($\Xi_i(v_h)$) to produce the noise-free MRI data (c);
 the MRI grid contains $N=9$ VPD.
 The Gaussian noise model with standard deviation equal to $\kappa=20\%$ of
 the peak velocity is added to the noise-free MRI data to produce the actual
 MRI data ($\bar{v}_i$) (d). Simulation-based imaging is
 used to reconstruct the velocity field from
 the noisy MRI data, which leads to the field 
 ($v_H(\,\cdot\,;\mu^\star)$) (e) and corresponding
 representation in the MR data space ($\Xi_i(v_H(\,\cdot\,;\mu^\star))$) (f).
 Colorbar: $\norm{v}$ [cm/s].}
 \label{fig:stenosis_main}
\end{figure}

\begin{figure}
 \centering
 \begin{tikzpicture}
\begin{groupplot} [
group style={group size = 3 by 2, horizontal sep = 0.5cm, vertical sep = 0.6cm},
title style={at={(current bounding box.north west)}, anchor=west}]
\nextgroupplot[axis equal image, axis lines=none, width=0.5\textwidth, ymax=0.1186, xmax=0.13, title=(a) Reference mesh, xmin=0, ymin=-0.079]
\addplot []
graphics [xmin=0,xmax=0.1248,ymin=-0.079,ymax=0.1186] { 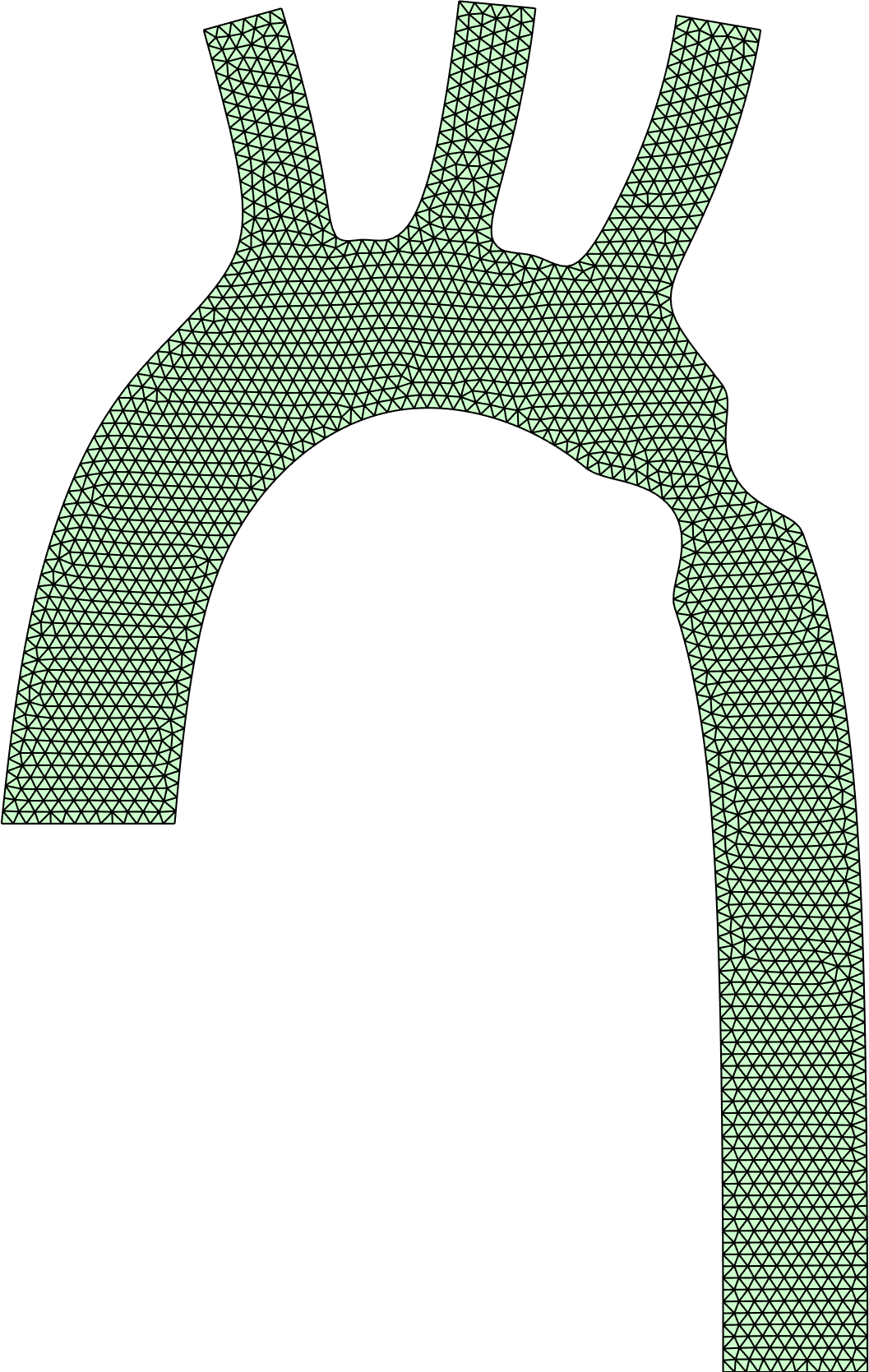};

\nextgroupplot[axis equal image, axis lines=none, width=0.5\textwidth, ymax=0.1186, xmax=0.13, title=(b) Reference velocity, xmin=0, ymin=-0.079]
\addplot []
graphics [xmin=0,xmax=0.1248,ymin=-0.079,ymax=0.1186] { 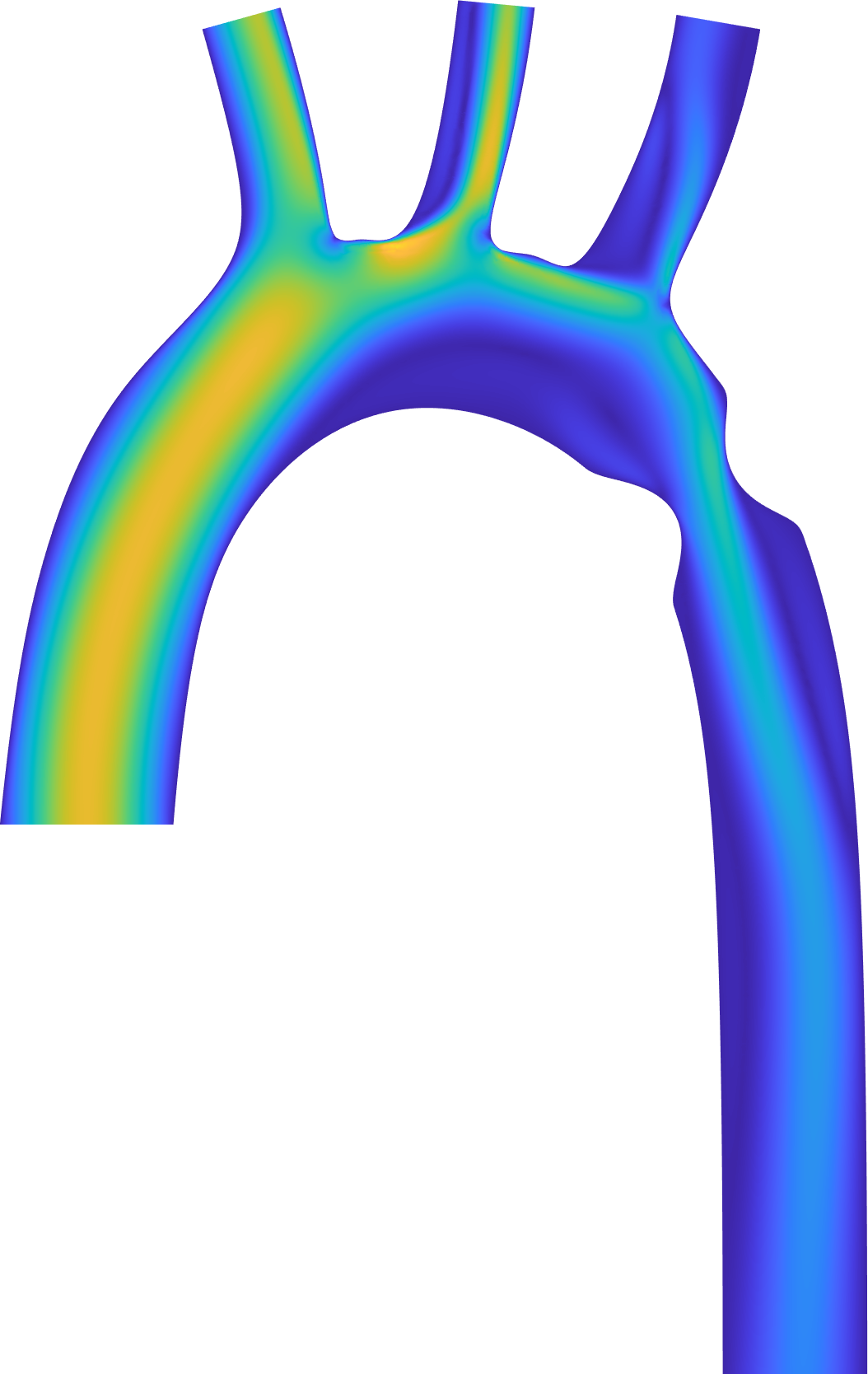};

\nextgroupplot[axis equal image, axis lines=none, width=0.5\textwidth, ymax=0.1186, xmax=0.13, title=(c) MR data (no noise), xmin=0, ymin=-0.079]
\addplot []
graphics [xmin=0,xmax=0.13,ymin=-0.079,ymax=0.116] { 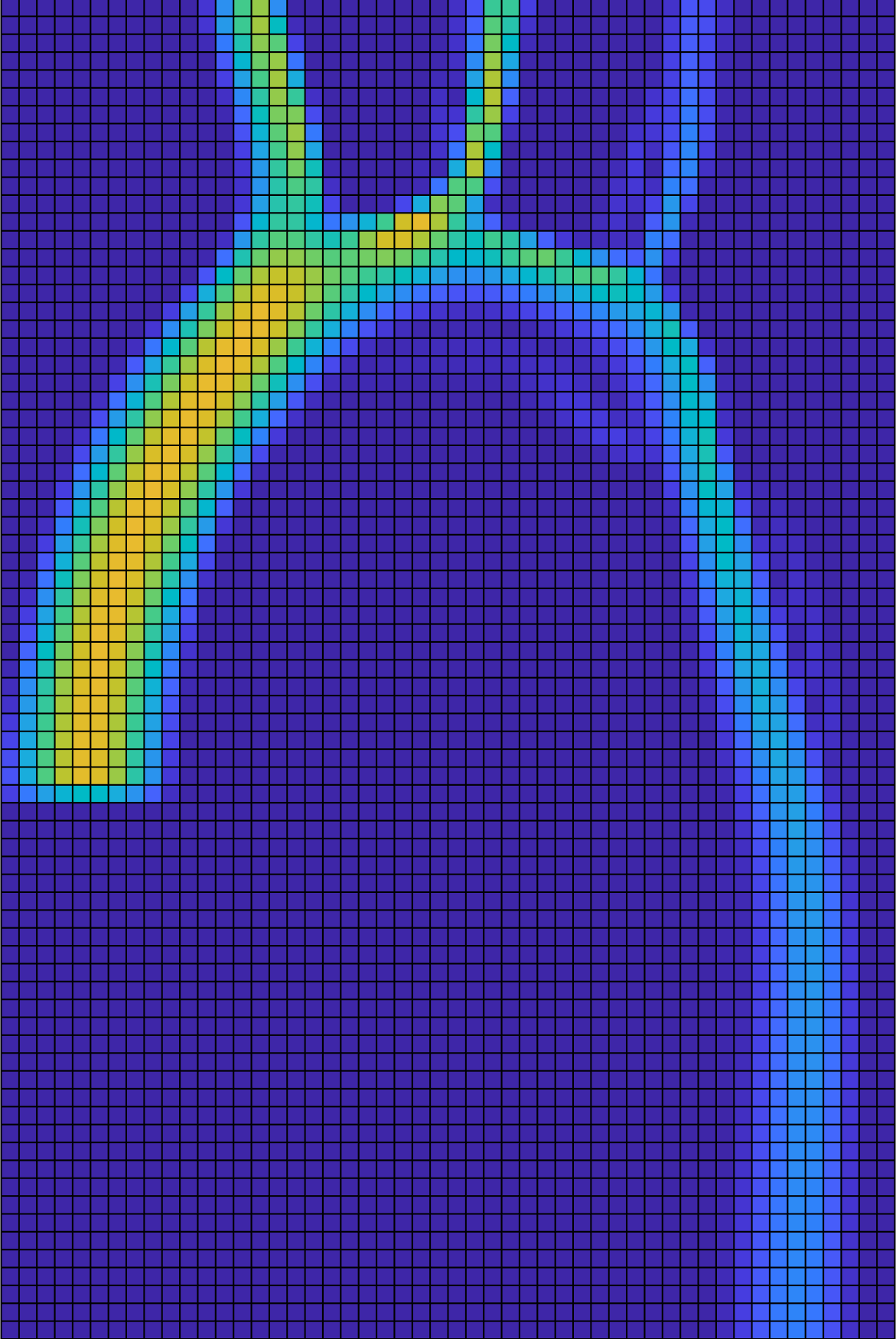};

\nextgroupplot[axis equal image, axis lines=none, width=0.5\textwidth, ymax=0.1186, xmax=0.13, title=(d) MR data (noise), xmin=0, ymin=-0.079]
\addplot []
graphics [xmin=0,xmax=0.13,ymin=-0.079,ymax=0.116] { 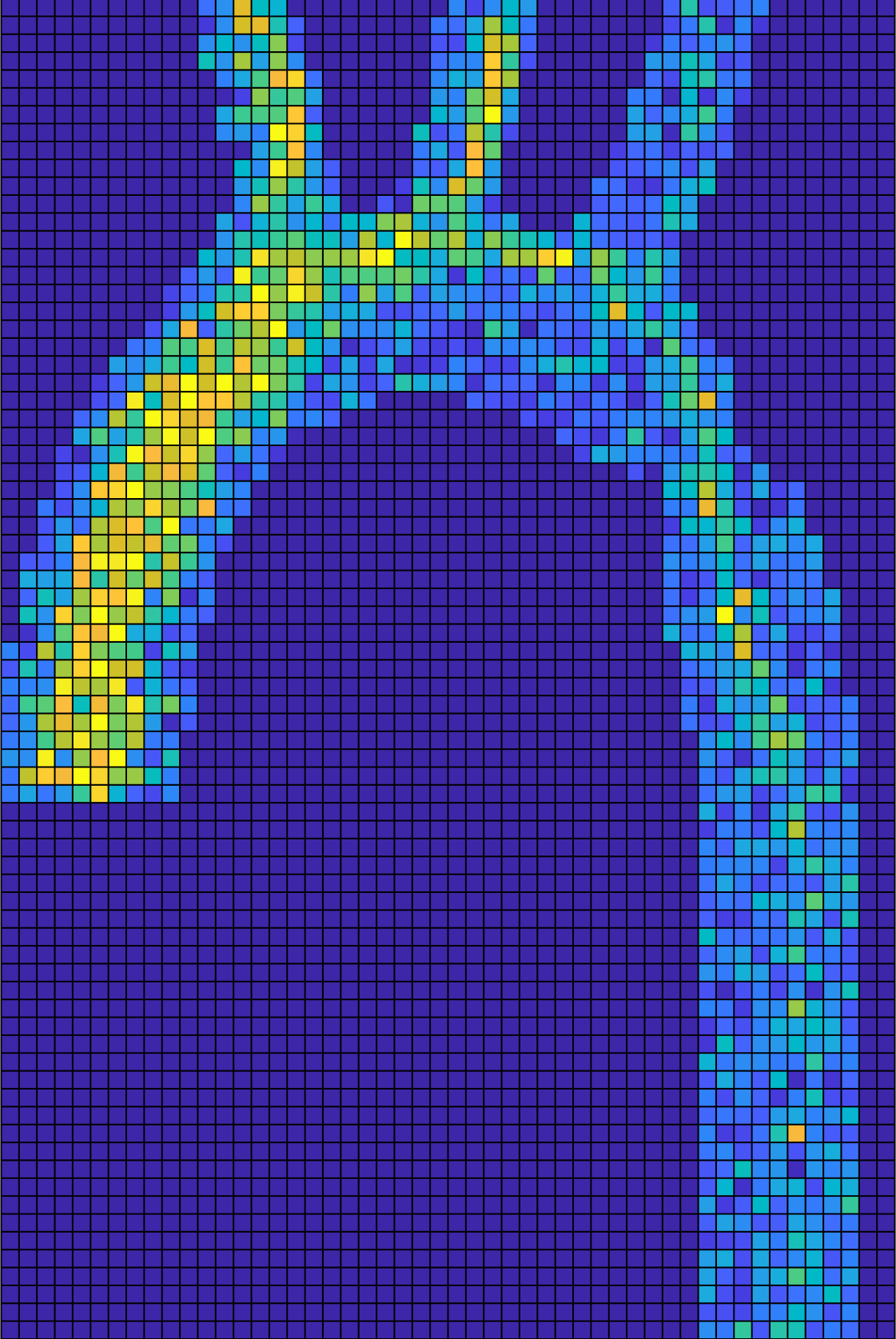};

\nextgroupplot[axis equal image, axis lines=none, width=0.5\textwidth, ymax=0.1186, xmax=0.13, title=(e) SBI reconstructed velocity, xmin=0, ymin=-0.079]
\addplot []
graphics [xmin=0,xmax=0.1248,ymin=-0.079,ymax=0.1186] { 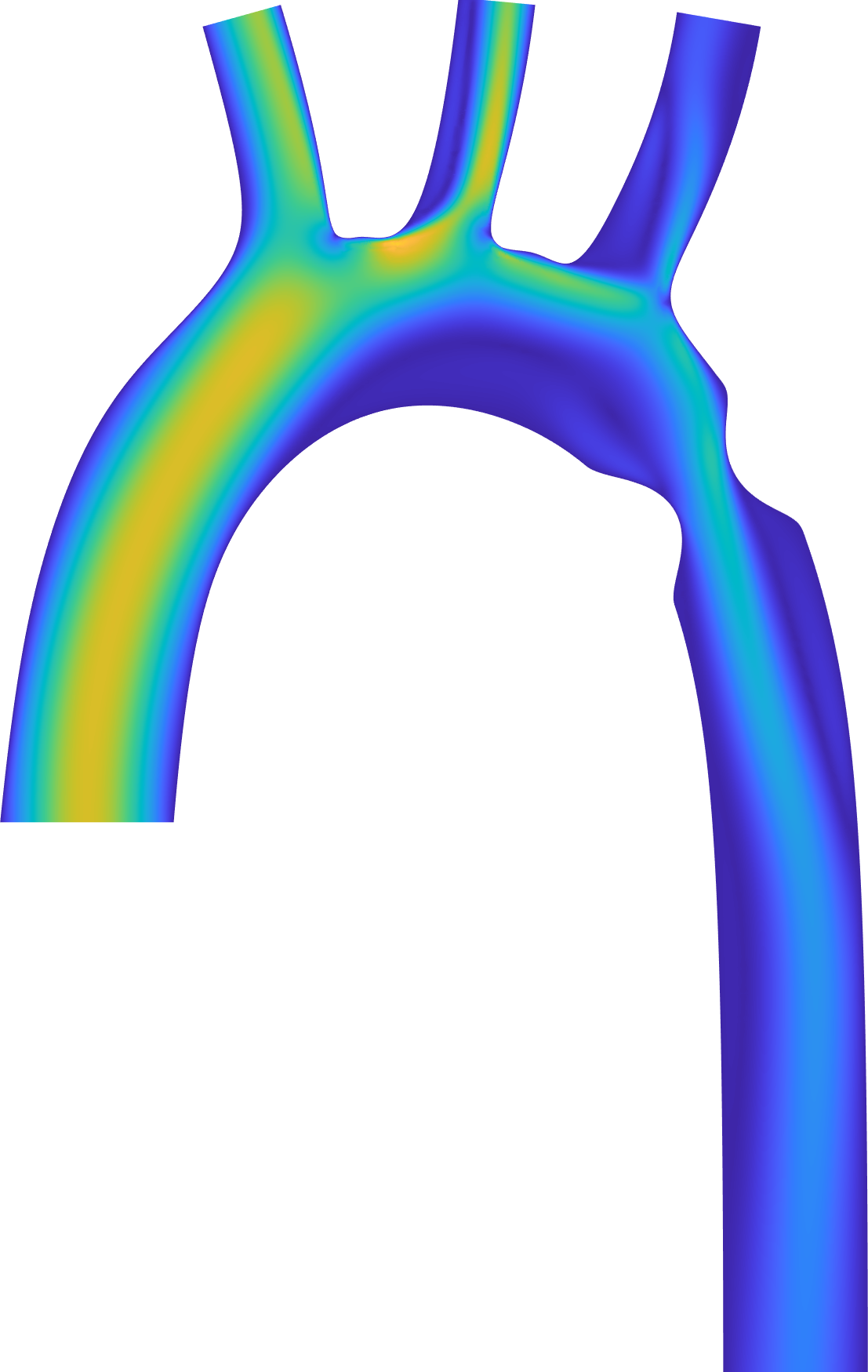};

\nextgroupplot[axis equal image, axis lines=none, width=0.5\textwidth, ymax=0.1186, xmax=0.13, title=(f) MR data from SBI, xmin=0, ymin=-0.079]
\addplot []
graphics [xmin=0,xmax=0.13,ymin=-0.079,ymax=0.116] { 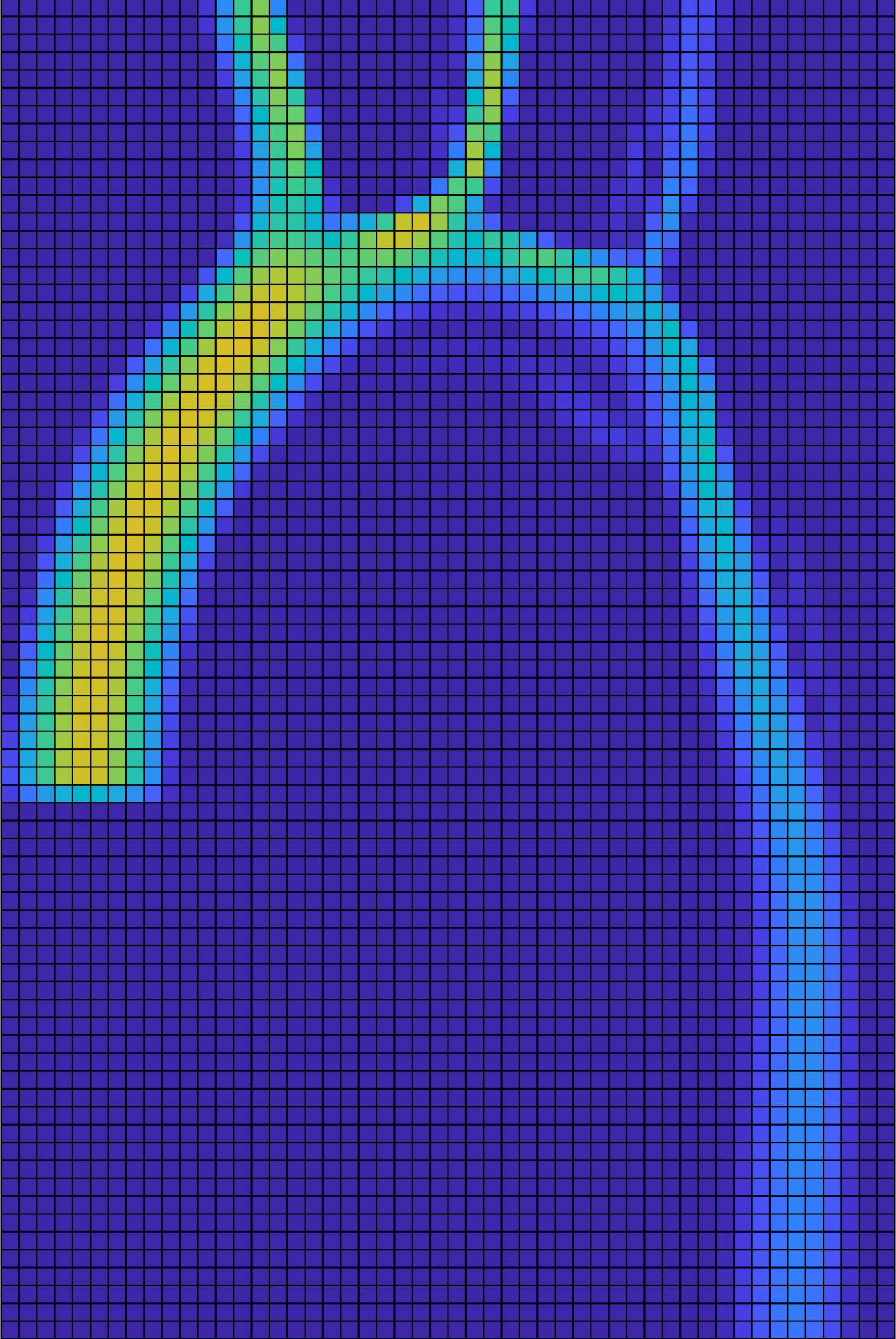};

\end{groupplot}\end{tikzpicture}
 \colorbarMatlabParula{0}{3.56}{7.15}{10.7}{14.3}
 \caption{Computational mesh (a) and corresponding velocity
 field (b) used to define the true flow ($v_h$) through the
 aorta ($\mathrm{Re}=1000$).
 The velocity field is mapped to the MR data space ($\Xi_i(v_h)$) to produce
 the noise-free MRI data (c); the MRI grid contains $N=10$ VPD.
 The Gaussian noise model with standard deviation equal to $\kappa=20\%$ of
 the peak velocity is added to the noise-free MRI data to produce the actual
 MRI data ($\bar{v}_i$) (d). Simulation-based imaging is
 used to reconstruct the velocity field from
 the noisy MRI data, which leads to the field
 ($v_H(\,\cdot\,;\mu^\star)$) (e) and corresponding
 representation in the MR data space ($\Xi_i(v_H(\,\cdot\,;\mu^\star))$) (f).
 Colorbar: $\norm{v}$ [cm/s].}
 \label{fig:aorta_main}
\end{figure}

\begin{figure}
 \centering
 \input{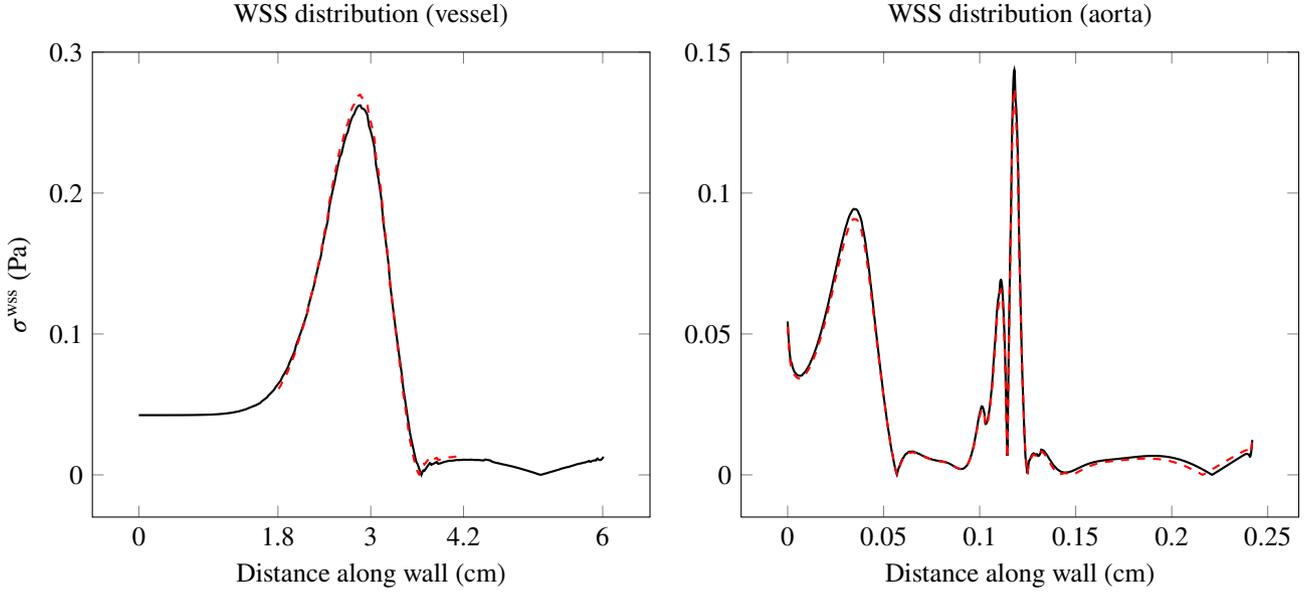}
 \caption{True wall shear stress distribution $\sigma_h^\mathrm{wss}$
          (\ref{line:wss_truth})
          and its reconstruction using simulation-based imaging
          $\sigma_H^\mathrm{wss}(\,\cdot\,,\mu^\star)$ from MRI grids ($N=9$ VPD
          for stenosis, $N=10$ VPD for aorta) with a noise level of
          $\kappa=20\%$ (\ref{line:wss_sbi}) along the intersection of
          the top wall of stenotic vessel (\textit{left}) and bottom wall
          of the aorta (\textit{right}) with the MRI domain,
          both at $\mathrm{Re}=1000$.}
 \label{fig:wss_truth_sbi}
\end{figure}

The Navier-Stokes equations are approximated using the finite element method
on an unstructured triangular mesh consisting of
$\mathcal{P}^3-\mathcal{P}^2$ Taylor-Hood elements. A linear mesh
(straight-sided triangles) is generated using DistMesh \cite{persson2004simple}
and the boundary edges are projected onto the exact geometry (and interior
nodes smoothed) for a high-order representation. Let $v_h$ denote the finite
element solution of the Navier-Stokes equations. This CFD solution defines our
reference or ``truth'' flow, e.g., corresponding to the \textit{in vivo} flow,
which is not available in practice, but essential to conduct thorough inquiries.
The reference or ``truth'' value for WSS is computed from (\ref{eqn:wss}) using
the finite element solution $v_h$; the pointwise velocity and necessary
derivatives are readily available from the finite element basis functions.
The reference WSS value computed from the reference finite element solution
$v_h$ is denoted $\sigma_h^\mathrm{wss}$. Figures~\ref{fig:stenosis_main}
and \ref{fig:aorta_main} show the computational mesh and corresponding
velocity field ($\mathrm{Re}=1000$) used to define the true flow for
the stenosis and aorta test cases, respectively. The corresponding WSS
($\sigma_h^\mathrm{wss}$) is shown in Figure~\ref{fig:wss_truth_sbi} for
both test cases. Synthetic MRI data is
extracted from the reference solution $v_h$ and perturbed with noise
using the approach in \cite{toger_zahr_4Dflow} (summarized in
Section~\ref{sec:methods:mri}).
In Section~\ref{sec:results}, we will use this synthetic MRI data to
reconstruct the WSS using standard MRI techniques
(Section~\ref{sec:methods:mri}) and SBI (Section~\ref{sec:methods:sbi})
to study the accuracy and sensitivity of each approach.

\subsection{Magnetic resonance imaging}
\label{sec:methods:mri}
Flow MRI scans extract velocities averaged over a Cartesian grid of voxels from an \textit{in vivo} flow. The resolution of the voxel grid determines both the resolution and noise of the velocity data \cite{edelstein1986intrinsic,portnoy2009information}. In blood flow imaging, the resolution is measured and reported in millimeters, but an important metric for flow accuracy is the number of voxels per diameter (VPD) \cite{wolf1993analysis}, which can range from 5-20 VPD for adults \cite{garcia2018distribution} and 3-5 VPD for infants \cite{akturk2018normal,cebral2009hemodynamics}. 

In the synthetic setting, we compute synthetic data consistent with the
approach in \cite{toger_zahr_4Dflow}, specialized to the case of steady flow,
i.e., a weighted integral of the true velocity field over a given voxel and
its neighbors. For simplicity, we assume the voxel grid is aligned with the
coordinate axes. We consider a grid consisting of $N_x$ voxels in the
$x_1$-direction and $N_y$ voxels in the $x_2$-direction and let
$\Delta x, \Delta y \in \Rbb_{>0}$ denote the spacing of the voxel grid in
the respective direction. We endow the $N=N_xN_y$ voxels with an ordering and
let $(X_i,Y_i)\in\Rbb^2$ denote the centroid of the $i$th voxel for
$i=1,\dots,N$. We leverage an abuse of notation to let $N$ denote the
resolution of the voxel grid, either stated in terms of the total number
of voxels ($N_xN_y$) or the number of VPD, depending on the context.
With this notation, the synthetic MR flow velocity data
associated with the $i$th voxel, $\bar{v}_i\in\Rbb^2$, is extracted from a
CFD simulation as
\begin{equation} \label{eqn:point_spread}
    \bar{v}_i \coloneqq \Xi_i(v_h) + \varphi_i,
\end{equation}
where $\varphi_i$ is a normally distributed random variable with mean $0$ and standard deviation proportional to the peak flow velocity, i.e., $\kappa\sup_{x\in\Omega} v_h(x)$ with noise level $\kappa\in\Rbb_{\geq 0}$, and
$\varphi_1,\dots,\varphi_N$ are independent and identically distributed
\cite{aja2013review,mcgibney1993unbiased}.
Following the approach in \cite{toger_zahr_4Dflow}, the point-spread 
function maps a continuous velocity field to the MR data space through a weighted average of the velocity field over a given voxel as
\begin{equation} \label{eqn:velmap}
    \Xi_i : u \mapsto \int_\Omega w_i u \, dv
\end{equation}
for $i=1,\dots,N$.
The weighting function for the $i$th voxel,
$\func{w_i}{\Omega}{\Rbb}$, is the normalized tensor product of a sinc
function with a smoothed box centered at $(X_i,Y_i)$, i.e.,
\begin{equation}
    w_i : x \mapsto
    \bar{w}_i(x)\left(\int_\Omega \bar{w}_i \, dv\right)^{-1}, \qquad
    \bar{w}_i : x \mapsto
    \Psi(x_1,X_i,\Delta x)
    \Psi(x_2,Y_i,\Delta y),
\end{equation}
where the component-wise, non-normalized weighting function, \func{\Psi}{\Rbb\times\Rbb\times\Rbb_{>0}}{\Rbb}, is defined as
\begin{equation}
    \Psi : (s, c, \Delta s) \mapsto
    \mathrm{sinc}\left(\frac{s-c}{\Delta s}\right) \chi(s, c, 4\Delta s).
\end{equation}
The sinc function is included to mimic the point-spread function of MRI
scanners that use a Fourier transform to map raw MRI data into flow velocities.
The one-dimensional smoothed box function, $\func{\chi}{\Rbb\times\Rbb\times\Rbb_{>0}}{\Rbb}$, is defined as
\begin{equation}
    \chi : (s, s_0, \omega) \mapsto
    \frac{1}{1+\mathrm{exp}(-(s-(s_0-\omega/2))/\gamma)}-
    \frac{1}{1+\mathrm{exp}(-(s-(s_0+\omega/2))/\gamma)}
\end{equation}
with center $s_0$, width $\omega$, and smoothness parameter $\gamma$.
The smoothed box function localizes the integrand in (\ref{eqn:velmap}) to the
center of a particular voxel, allowing for some overlap between voxels, and
ensures the integrand is sufficiently smooth for the integral to be well-approximated using numerical quadrature.
Following the work in \cite{toger_zahr_4Dflow}, we take the smoothness parameter to be
proportional to the voxel spacing, $\gamma = 0.1\min\{\Delta x, \Delta y\}$.
In this work, (\ref{eqn:point_spread}) is used to define the synthetic MRI data and
define the SBI cost function (Section~\ref{sec:methods:sbi}); a complete
description of the approach can be found in \cite{toger_zahr_4Dflow}.
Figures~\ref{fig:stenosis_main} and \ref{fig:aorta_main} show the
synthetic MRI data ($N=9$ VPD for stenosis, $N=10$ VPD for aorta)
extracted from the true flow with ($\kappa=20\%$) and without noise
for the stenosis and aorta test cases, respectively.

We use standard methods from the MRI literature to compute the WSS directly
from the MRI data,
$\sigma_N^\mathrm{wss}(x)$, at any point $x\in\partial\Omega_\text{w}$
\cite{quant_wss_numerical,osinnski1995determination}.
First, we use the raw voxel velocity data $v_h$ to compute a bilinear
flow field reconstruction; we call this bilinear flow field $v_N$. Then,
assuming no flow penetration through the walls, the tangential component of
the surface traction in (\ref{eqn:wss}) can be written as
\begin{equation} \label{eqn:wss_tract_simple}
 \tau = \mu (BB^T) \nabla v \cdot n,
\end{equation}
where the columns of $B(x)\in\Rbb^{d\times (d-1)}$ form a basis of the
tangent space at $x\in\Omega_\mathrm{w}$ (orthogonal complement of $n(x)$);
see derivation in Appendix~\ref{app:mriwss}. The normal gradient is computed
using an approach similar to that in
\cite{quant_wss_numerical,osinnski1995determination},
i.e., construct a quadratic approximation of the velocity field in the normal
direction from the bilinear flow field and assume the velocity is zero at $x$,
because it performed favorably relative to several alternatives in a
comparative study \cite{quant_wss_numerical}. The quadratic velocity along
the outward normal $n$ originating at the point $x$ takes the form
\begin{equation}
 \func{\tilde{v}_N(\,\cdot\,; x, n)}{\Rbb}{\Rbb^d}, \qquad
 \tilde{v}_N(\,\cdot\,; x, n) : r \mapsto
 v_N(x-\delta n)\psi_1(r) + v_N(x-2\delta n) \psi_2(r)
\end{equation}
where $\{\psi_0,\psi_1,\psi_2\}$ are one-dimensional quadratic Lagrangian polynomials associated with the nodes $\{0, \delta, 2\delta\}$, and $\delta\in\Rbb_{>0}$ is the increment used to determine the points along the normal at which to fit the quadratic function to the bilinear flow field. In this work, we choose the increment to be proportional to the voxel spacing, $\delta=1.2 \min\{\Delta x, \Delta y\,, 0.06\}$ (cm). Finally, the gradient of the velocity in the normal direction is approximated as $\tilde{v}_N'(0)$ and the tangential component of the surface traction ($\tau_N$) and wall shear stress ($\sigma_N^\mathrm{wss}$) at $x\in\Omega$ are computed as
\begin{equation}
 \tau_N(x) = \mu(B(x)B(x)^T) \tilde{v}_N'(0; x, n(x)), \qquad
 \sigma_N^{\mathrm{wss}}(x) = \norm{\tau_N(x)}.
\end{equation}

This approach makes two unrealistic assumptions, namely, that the point on the boundary $x$ and corresponding normal $n(x)$ are known exactly. In practice, these geometric properties must be approximated from scanned images, which introduces additional error that was quantified in \cite{quant_wss_numerical}. Since the position on the wall and the corresponding normal are known in the SBI setting, we use this information to maintain fairness in the comparison between MRI postprocessing and SBI.

\subsection{Simulation-based imaging}
\label{sec:methods:sbi}
Simulation-based imaging aims to reconstruct a high-fidelity \textit{in vivo} flow image from a CFD simulation that has been certified with MRI flow measurements. It optimally fits a CFD simulation to MRI flow data that can be noisy, sparse, and low-resolution by modifying the boundary conditions, material properties, and the initial condition. In this work, we adjust the inflow boundary conditions to fit the CFD simulation to the MRI data. That is, we consider the Navier-Stokes equation in (\ref{eqn:ins}) subject to the following boundary conditions
\begin{equation} \label{eqn:bc1}
    v = 0~~\text{on}~~\partial\Omega_\text{w}, \quad
    v = \hat{v}(\,\cdot\,;\mu)~~\text{on}~~\partial\Omega_\text{in}, \quad
    \sigma \cdot n = 0~~\text{on}~~\partial\Omega_\text{out}
\end{equation}
where $\func{\hat{v}}{\Rbb^d\times \Dcal}{\Rbb^d}$ is the parametrized inflow function, $\mu\in\Dcal$ is a vector of parameters, and $\Dcal\subset\Rbb^d$ is the admissible parameter space. In this work, we take the inflow velocity to be parallel to the normal of the
inflow boundary surface following \cite{toger_zahr_4Dflow}, which leads to
\begin{equation}
    \hat{v} : (x;\mu) \mapsto \frac{(B_0 - x_2)(B_0 + x_2)}{B_0^2} (\mu, 0)
\end{equation}
for the vessel case study,
i.e., a parabolic profile for the $x_1$ velocity that is zero at the wall ($x_2=\pm B_0$) with $\mu$ defining the peak of the parabola. A similar parabolic parametrization of the normal-directed inflow velocity is used for the aorta.

The CFD simulation underlying SBI discretizes the Navier-Stokes equation in
(\ref{eqn:ins}) with the parametrized boundary
conditions in (\ref{eqn:bc1}) using the finite element method as described in Section~\ref{sec:methods:synthetic}; the corresponding velocity field is denoted $v_H(x;\mu)$.
In the \textit{in vivo} setting, the geometry of the flow domain is obtained by segmenting an angiogram scan from which a mesh is generated using standard tools \cite{kang2012heart,prakosa2014methodology,stanescu2010investigation}.
In this study, we directly generate a high-order mesh of the two geometries considered (Section~\ref{sec:methods:synthetic}). 

The parameterized inflow boundary conditions are determined by optimally fitting the parameterized CFD solution to the MRI data 
\begin{equation} \label{eqn:opt_prob}
        \mu^\star = \argmin_{\mu\in\Dcal} ~I(\mu), \qquad I : \mu \mapsto \sum_{i=1}^N \frac{\alpha_i}{2} \norm{\Xi_i(v_H(\,\cdot\,;\mu)) - \bar{v}_i}_2^2,
\end{equation}
where $\alpha_i = 1$ if the $i$th voxel at least partially within the domain
and zero otherwise and $\func{I}{\Dcal}{\Rbb}$ is the cost function that
measures the misfit between the MRI data and its prediction from the CFD
simulation. The optimization problem in
(\ref{eqn:opt_prob}) is solved using a quasi-Newton method globalized with a
line search \cite{nocedal2006numerical} and gradients of the objective
function are computed efficiently using the adjoint method.
From the solution of the SBI optimization problem
($\mu^\star$), the reconstructed flow field is the CFD simulation at the
optimal parameter configuration, i.e., $v_H(\,\cdot\,;\mu^\star)$. From the
SBI reconstructed flow, we calculate the corresponding WSS, denoted
$\sigma_H^\mathrm{wss}(\,\cdot\,;\mu^\star)$, using (\ref{eqn:wss}) with
the SBI state $v_H(\,\cdot\,; \mu^\star)$. The necessary derivatives of the
flow solution required in (\ref{eqn:wss}) are readily available from the
finite element basis functions.

Figures~\ref{fig:stenosis_main} and \ref{fig:aorta_main} show the SBI
reconstruction of the flow field and its representation in the MR
data space for the stenosis and aorta test cases, respectively.
Additionally, Figure~\ref{fig:wss_truth_sbi} shows the WSS
reconstruction using SBI for both test cases.
From these figures, it is clear that even with a moderately refined
MRI voxel mesh ($N=9$ VPD for stenosis, $N=10$ VPD for aorta)
and high noise level ($\kappa=20\%$), the SBI reconstruction provides
a good approximation to the true velocity field and WSS distribution.
We will study its sensitivity to these parameters in Section~\ref{sec:results}.

The same point-spread function ($\Xi$) used to define the MRI data is used to sample the CFD solution for matching to the MR flow data in the objective function. This assumes the numerical point-spread function (\ref{eqn:point_spread}) will exactly reproduce the point-spread function of the MRI scanner, which is not true in practice, which may introduce additional modeling error. However, in this work, we do not consider the sensitivity of SBI to the point-spread function, instead we focus on its performance with respect to Reynolds number, MRI resolution, and noise.

\section{Results}
\label{sec:results}
In this section, we study the performance of MRI- and SBI-based wall shear
stress reconstructions as a function of MRI noise level, resolution of the MRI
voxel grid, Reynolds number, and resolution of the CFD mesh used for SBI
reconstruction. 
Comprehensive studies for the ideal stenosis geometry are
used to draw conclusions regarding the performance of MRI and SBI wall shear
stress reconstruction; these conclusions are then verified for the more complex
setting of flow through an idealized aorta using targeted studies.

In each of these studies, we will quantitatively compare the WSS distributions,
i.e., $\sigma_H^\mathrm{wss}(\,\cdot\,;\mu^\star)$ (SBI) and
$\sigma_N^\mathrm{wss}$ (MRI) to $\sigma_h^\mathrm{wss}$ (true WSS),
along the curve $\Gamma\subset\partial\Omega_\text{w}$, where $\Gamma$
is the intersection of the top (bottom) wall of the vessel (aorta) with
the limits of the MRI domain (which does not span the entire domain for the
stenotic vessel case). The error will be quantified using the relative $L^2$
norm, reported as percentages,
\begin{equation} \label{eqn:comp_met}
    e_{\mathrm{SBI}} = \frac{\sqrt{\int_\Gamma |\sigma_h^\mathrm{wss}(x)-\sigma_H^\mathrm{wss}(x;\mu^\star)|^2\,dS}}{\sqrt{\int_\Gamma |\sigma_h^\mathrm{wss}(x)|^2\,dS}}, \qquad
    e_{\mathrm{MRI}} = \frac{\sqrt{\int_\Gamma |\sigma^\mathrm{wss}_h(x)-\sigma^\mathrm{wss}_N(x)|^2\,dS}}{\sqrt{\int_\Gamma |\sigma^\mathrm{wss}_h(x)|^2\,dS}},
\end{equation}
where $e_{\mathrm{SBI}}, e_{\mathrm{MRI}}\in\Rbb_{\geq 0}$ are the SBI and MRI WSS distribution errors,
respectively, and the integrals are computed using Gaussian quadrature.


\subsection{Impact of the CFD resolution}
\begin{figure}
 \centering
 \begin{tikzpicture}
\begin{groupplot} [
group style={group size = 2 by 3, horizontal sep = 0.4cm, vertical sep = 0.2cm}]
\nextgroupplot[axis equal image, axis lines=none, width=0.58\textwidth, ymax=0.05, xmax=1, xmin=0, ymin=-0.05]
\addplot []
graphics [xmin=0,xmax=1,ymin=-0.05,ymax=0.05] { 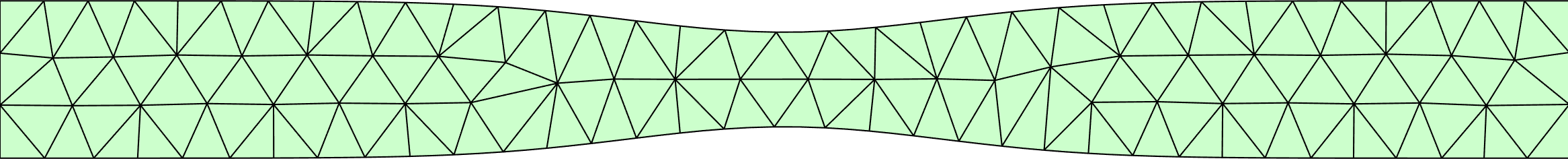};

\nextgroupplot[axis equal image, axis lines=none, width=0.58\textwidth, ymax=0.05, xmax=1, xmin=0, ymin=-0.05]
\addplot []
graphics [xmin=0,xmax=1,ymin=-0.05,ymax=0.05] { 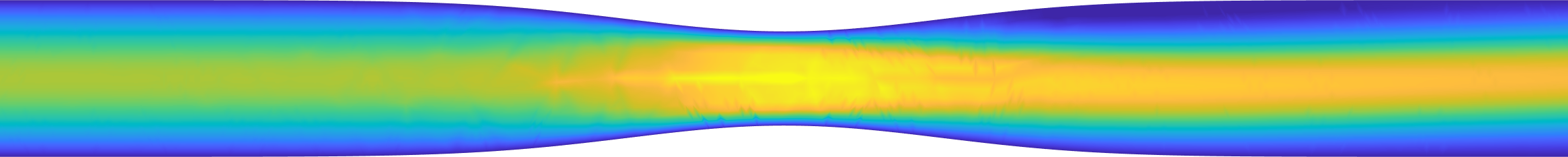};

\nextgroupplot[axis equal image, axis lines=none, width=0.58\textwidth, ymax=0.05, xmax=1, xmin=0, ymin=-0.05]
\addplot []
graphics [xmin=0,xmax=1,ymin=-0.05,ymax=0.05] { 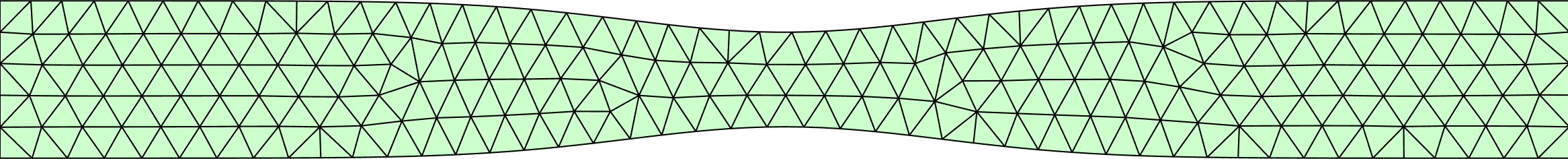};

\nextgroupplot[axis equal image, axis lines=none, width=0.58\textwidth, ymax=0.05, xmax=1, xmin=0, ymin=-0.05]
\addplot []
graphics [xmin=0,xmax=1,ymin=-0.05,ymax=0.05] { 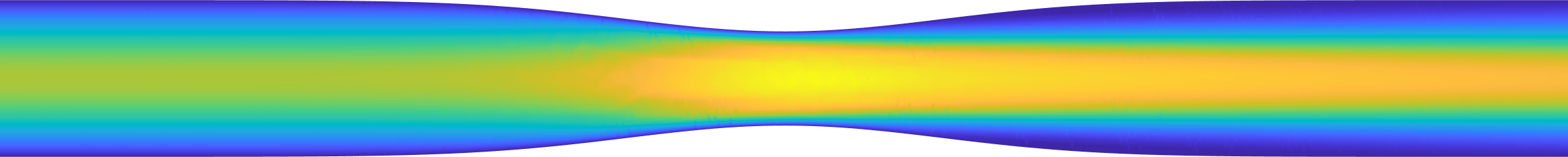};

\nextgroupplot[axis equal image, axis lines=none, width=0.58\textwidth, ymax=0.05, xmax=1, xmin=0, ymin=-0.05]
\addplot []
graphics [xmin=0,xmax=1,ymin=-0.05,ymax=0.05] { 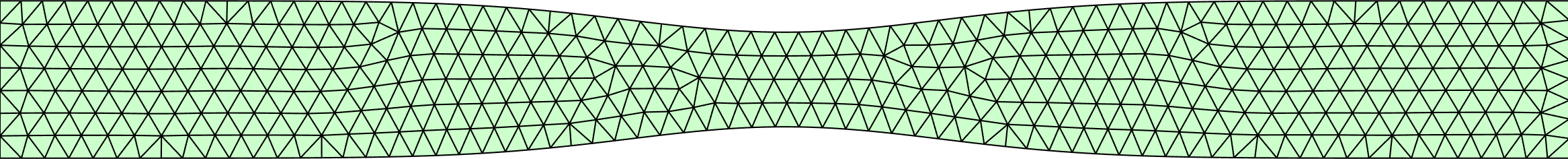};

\nextgroupplot[axis equal image, axis lines=none, width=0.58\textwidth, ymax=0.05, xmax=1, xmin=0, ymin=-0.05]
\addplot []
graphics [xmin=0,xmax=1,ymin=-0.05,ymax=0.05] { 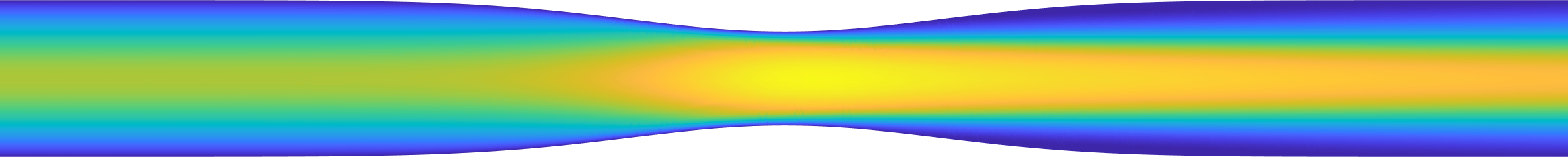};

\end{groupplot}\end{tikzpicture}
 \caption{Meshes (\textit{left}) of stenotic vessel used to study the
 sensitivity of SBI to resolution of the CFD mesh and the corresponding
 velocity field (\textit{right}). Number of elements in mesh:
 $N_e=368$ (\textit{top row}), $N_e=766$ (\textit{middle row}), and
 $N_e=1590$ (\textit{bottom row}). Colorbar in Figure~\ref{fig:stenosis_main}.}
 \label{fig:cfd_meshes}
\end{figure}

\ \label{sec:CFDres}
The CFD mesh resolution used for SBI impacts the computational cost
of the method so we study the sensitivity of the SBI wall
shear stress reconstruction to the resolution of the mesh. For this
study, we vary the noise $\kappa\in\{0, 5, 10, 15, 20\}$ ($\%$), Reynolds number
$\mathrm{Re}\in\{100,500,1000\}$, MRI grid resolution $N\in\{3, 9, 15, 28\}$
(VPD), and consider the three meshes of increasing resolution shown in
Figure~\ref{fig:cfd_meshes} ($N_e = 368$, $766$, $1590$ elements, respectively).
The finest mesh ($N_e=1590$) is used to compute the reference (``true'') flow and smoothly
resolves all flow features, whereas the two coarser meshes
($N_e=368$ and $N_e=766$) lead to numerical artifacts in the
stenosis, but contain two and four times fewer elements, respectively,
and the corresponding CFD simulations require a fraction of the compute time. 

\begin{figure}
 \centering
 \input{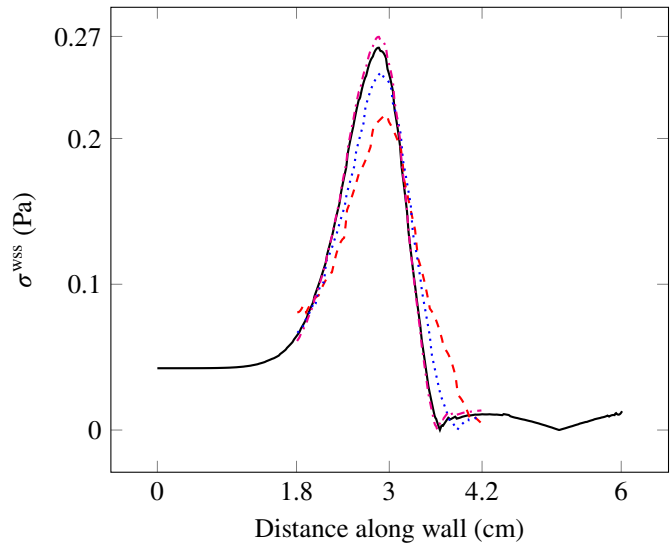}
 \caption{Wall shear stress distribution over intersection of top wall of
  stenotic vessel with MRI domain ($\Gamma$) using SBI with different
  mesh resolutions at Reynolds number $\mathrm{Re}=1000$, noise
  $\kappa=20\%$, and MRI resolution $N=9$ VPD (scenario in
  Figure~\ref{fig:stenosis_main}). Legend: WSS distribution
  from true flow (\ref{line:wss_truth}) and SBI WSS distribution using the
  mesh with $N_e=368$ elements (\ref{line:wss_sbi_nref0}), $N_e=766$ elements
  (\ref{line:wss_sbi_nref1}), and $N_e=1590$ elements
  (\ref{line:wss_sbi_nref2}); see Figure~\ref{fig:cfd_meshes} for meshes.}
 \label{fig:cfd_trials_wss}
\end{figure}

First, we observe that the WSS reconstruction error using SBI decreases
as the CFD mesh is refined (Figure~\ref{fig:sten_err_wss_cfd}), which
holds for almost all Reynolds numbers, noise levels, and MRI voxel grids considered.
The noise in the MRI data tends to degrade the accuracy of the WSS
reconstruction; however, its influence diminishes as the MRI voxel
grid is refined, i.e., there is very little difference in the $\kappa=0\%$
and $\kappa=20\%$ WSS reconstruction with $N=28$ VPD MRI grid, whereas there
is a substantial difference with $N=3$ VPD. We
note that for the finest mesh ($N_e=1590$), the WSS reconstruction has a
relatively weak dependence on the Reynolds number of the flow, whereas
there is a significant dependence on Reynolds number for the coarser meshes,
e.g., for the coarse mesh with $N=3$ VPD and $\kappa=10\%$, the WSS error is
around $5\%$ for $\mathrm{Re}=100$ and close to $15\%$ for $\mathrm{Re}=1000$.
Even in the cases with large error, the optimization problem underlying SBI
drives the CFD simulation to a configuration that agrees reasonably well
with the true WSS distribution given the limitations of the discretization
(Figure~\ref{fig:cfd_trials_wss}).

Together these observations imply that an underresolved mesh can be
used for the SBI reconstruction at the cost of increasing the
sensitivity of the WSS reconstruction to Reynolds number and
the MRI voxel grid. This is significant because it
means the high computational cost of SBI can be reduced by using
a coarse mesh without significant loss in WSS reconstruction accuracy,
provided a sufficiently high-resolution MRI data set is used. Alternatively,
it implies that a limited amount of MRI data can be used (e.g., from fast
patient scans) provided the CFD mesh used for SBI reconstruction is relatively
fine, which indicates an inherent trade-off between scan and reconstruction time
when using SBI. In the remainder, we will fix the CFD mesh resolution at $N_e=1590$
and focus on the impact of the remaining parameters (noise, Reynolds number, MRI
resolution).


\begin{figure}
 \centering
 \input{_py/error_fcnof_cfd_newdat.tikz}
 \caption{WSS reconstruction error using SBI as a function of CFD mesh resolution for various MRI grids (\textit{columns}) and noise levels (\textit{lines}) at Reynolds number $100$ (\textit{top row}), $500$ (\textit{middle row}), and $1000$ (\textit{bottom row}).
 Legend: $\kappa=0\%$ (\ref{line:err_sbi_cfd0}),
         $\kappa=5\%$ (\ref{line:err_sbi_cfd1}),
         $\kappa=10\%$ (\ref{line:err_sbi_cfd2}),
         $\kappa=15\%$ (\ref{line:err_sbi_cfd3}),
         $\kappa=20\%$ (\ref{line:err_sbi_cfd4}).} 
 \label{fig:sten_err_wss_cfd}
\end{figure}

\subsection{Impact of noise, Reynolds number, and MRI resolution}
In this section, we study the coupled effect that noise, Reynolds number,
and the MRI voxel grid resolution have on the accuracy to which WSS is
reconstructed using SBI (Section~\ref{sec:methods:sbi}) and standard MRI
reconstruction (Section~\ref{sec:methods:mri}).
We vary the Reynolds number $\mathrm{Re}\in\{100,500,1000\}$ because it
is known to strongly impact the accuracy of WSS reconstructions
based solely on MRI \cite{quant_wss_numerical}.
We vary the noise $\kappa\in\{0, 5, 10, 15, 20\}$ ($\%$)
and MRI voxel grid $N\in\{3, 9, 15, 28\}$ (VPD) because these incorporate
extreme (best-case and worst-case) scenarios seen in practice;
usually, for infant patients, noise levels are $3\%-10\%$ and
MRI resolution is $3-5$ VPD \cite{akturk2018normal,cebral2009hemodynamics}.
We use the CFD mesh with $N_e=1590$ elements for SBI reconstruction.

\begin{figure}
 \centering
 \begin{tikzpicture}
\begin{groupplot} [
width=0.3\textwidth,
group style={group size = 4 by 2, horizontal sep = 0.3cm, vertical sep = 0.4cm}]
\nextgroupplot[xtick={3, 9, 15, 28}, ymin=0, ymax=10, ylabel={$e_{\mathrm{SBI}}$ (\%)}, title={$\kappa=0\%$}, xticklabels={,,}]
\addplot [black, solid, mark options={solid, thin}, mark=*, mark size=1.5]
coordinates {
( 3.00000000e+00,  3.55486846e-02)
( 9.00000000e+00,  2.11432050e-02)
( 1.50000000e+01,  3.96954920e-05)
( 2.80000000e+01,  1.92496800e-04)};\label{line:err_sbi_vpd0}

\addplot [blue, solid, mark options={solid, thin}, mark=square*, mark size=1.5]
coordinates {
( 3.00000000e+00,  5.72782673e-04)
( 9.00000000e+00,  1.41276278e-05)
( 1.50000000e+01,  1.64484550e-05)
( 2.80000000e+01,  2.34543308e-05)};\label{line:err_sbi_vpd1}

\addplot [red, solid, mark options={solid, thin}, mark=triangle*, mark size=1.5]
coordinates {
( 3.00000000e+00,  8.09942723e-04)
( 9.00000000e+00,  1.37969706e-05)
( 1.50000000e+01,  2.06027349e-05)
( 2.80000000e+01,  5.83797456e-06)};\label{line:err_sbi_vpd2}

\nextgroupplot[yticklabels={,,}, xtick={3, 9, 15, 28}, ymax=10, title={$\kappa=5\%$}, ymin=0, xticklabels={,,}]
\addplot [black, solid, mark options={solid, thin}, mark=*, mark size=1.5, forget plot]
coordinates {
( 3.00000000e+00,  2.22485149e+00)
( 9.00000000e+00,  7.19548739e-01)
( 1.50000000e+01,  3.30628594e-01)
( 2.80000000e+01,  6.46523112e-02)};

\addplot [blue, solid, mark options={solid, thin}, mark=square*, mark size=1.5, forget plot]
coordinates {
( 3.00000000e+00,  2.35051306e+00)
( 9.00000000e+00,  6.91964759e-01)
( 1.50000000e+01,  3.55280957e-01)
( 2.80000000e+01,  4.41273842e-02)};

\addplot [red, solid, mark options={solid, thin}, mark=triangle*, mark size=1.5, forget plot]
coordinates {
( 3.00000000e+00,  2.36878840e+00)
( 9.00000000e+00,  6.90204685e-01)
( 1.50000000e+01,  3.69171269e-01)
( 2.80000000e+01,  4.26644615e-02)};

\nextgroupplot[yticklabels={,,}, xtick={3, 9, 15, 28}, ymax=10, title={$\kappa=10\%$}, ymin=0, xticklabels={,,}]
\addplot [black, solid, mark options={solid, thin}, mark=*, mark size=1.5, forget plot]
coordinates {
( 3.00000000e+00,  4.47509452e+00)
( 9.00000000e+00,  1.41869455e+00)
( 1.50000000e+01,  6.61457667e-01)
( 2.80000000e+01,  1.29114460e-01)};

\addplot [blue, solid, mark options={solid, thin}, mark=square*, mark size=1.5, forget plot]
coordinates {
( 3.00000000e+00,  4.68653267e+00)
( 9.00000000e+00,  1.38500199e+00)
( 1.50000000e+01,  7.11002174e-01)
( 2.80000000e+01,  8.82311725e-02)};

\addplot [red, solid, mark options={solid, thin}, mark=triangle*, mark size=1.5, forget plot]
coordinates {
( 3.00000000e+00,  4.72260732e+00)
( 9.00000000e+00,  1.38190010e+00)
( 1.50000000e+01,  7.38778108e-01)
( 2.80000000e+01,  8.53286940e-02)};

\nextgroupplot[yticklabels={,,}, xtick={3, 9, 15, 28}, ymax=10, title={$\kappa=20\%$}, ymin=0, xticklabels={,,}]
\addplot [black, solid, mark options={solid, thin}, mark=*, mark size=1.5, forget plot]
coordinates {
( 3.00000000e+00,  8.94502468e+00)
( 9.00000000e+00,  2.81972424e+00)
( 1.50000000e+01,  1.32360011e+00)
( 2.80000000e+01,  2.58045694e-01)};

\addplot [blue, solid, mark options={solid, thin}, mark=square*, mark size=1.5, forget plot]
coordinates {
( 3.00000000e+00,  9.31281610e+00)
( 9.00000000e+00,  2.77433231e+00)
( 1.50000000e+01,  1.42371425e+00)
( 2.80000000e+01,  1.76438350e-01)};

\addplot [red, solid, mark options={solid, thin}, mark=triangle*, mark size=1.5, forget plot]
coordinates {
( 3.00000000e+00,  9.38217183e+00)
( 9.00000000e+00,  2.76818846e+00)
( 1.50000000e+01,  1.47935571e+00)
( 2.80000000e+01,  1.70638973e-01)};

\nextgroupplot[xtick={3, 9, 15, 28}, ymin=0, xlabel={$N$ (VPD)}, ymax=100, ylabel={$e_{\mathrm{MRI}}$ (\%)}]
\addplot [black, solid, mark options={solid, thin}, mark=*, mark size=1.5, forget plot]
coordinates {
( 3.00000000e+00,  2.65715241e+01)
( 9.00000000e+00,  1.78628927e+01)
( 1.50000000e+01,  1.07892147e+01)
( 2.80000000e+01,  5.05153628e+00)};

\addplot [blue, solid, mark options={solid, thin}, mark=square*, mark size=1.5, forget plot]
coordinates {
( 3.00000000e+00,  4.91077861e+01)
( 9.00000000e+00,  3.95811298e+01)
( 1.50000000e+01,  2.77682548e+01)
( 2.80000000e+01,  1.38287456e+01)};

\addplot [red, solid, mark options={solid, thin}, mark=triangle*, mark size=1.5, forget plot]
coordinates {
( 3.00000000e+00,  5.85280535e+01)
( 9.00000000e+00,  5.03418769e+01)
( 1.50000000e+01,  3.84932763e+01)
( 2.80000000e+01,  2.07782254e+01)};

\nextgroupplot[yticklabels={,,}, xtick={3, 9, 15, 28}, ymin=0, xlabel={$N$ (VPD)}, ymax=100]
\addplot [black, solid, mark options={solid, thin}, mark=*, mark size=1.5, forget plot]
coordinates {
( 3.00000000e+00,  2.78437783e+01)
( 9.00000000e+00,  2.38641502e+01)
( 1.50000000e+01,  1.99183403e+01)
( 2.80000000e+01,  2.73248578e+01)};

\addplot [blue, solid, mark options={solid, thin}, mark=square*, mark size=1.5, forget plot]
coordinates {
( 3.00000000e+00,  4.70431007e+01)
( 9.00000000e+00,  4.13301593e+01)
( 1.50000000e+01,  3.08691039e+01)
( 2.80000000e+01,  2.31931672e+01)};

\addplot [red, solid, mark options={solid, thin}, mark=triangle*, mark size=1.5, forget plot]
coordinates {
( 3.00000000e+00,  5.64441610e+01)
( 9.00000000e+00,  5.11903870e+01)
( 1.50000000e+01,  4.03305203e+01)
( 2.80000000e+01,  2.58547775e+01)};

\nextgroupplot[yticklabels={,,}, xtick={3, 9, 15, 28}, ymin=0, xlabel={$N$ (VPD)}, ymax=100]
\addplot [black, solid, mark options={solid, thin}, mark=*, mark size=1.5, forget plot]
coordinates {
( 3.00000000e+00,  3.32812207e+01)
( 9.00000000e+00,  3.55183830e+01)
( 1.50000000e+01,  3.41795120e+01)
( 2.80000000e+01,  4.76494314e+01)};

\addplot [blue, solid, mark options={solid, thin}, mark=square*, mark size=1.5, forget plot]
coordinates {
( 3.00000000e+00,  4.61364859e+01)
( 9.00000000e+00,  4.51786292e+01)
( 1.50000000e+01,  3.66709606e+01)
( 2.80000000e+01,  3.66924359e+01)};

\addplot [red, solid, mark options={solid, thin}, mark=triangle*, mark size=1.5, forget plot]
coordinates {
( 3.00000000e+00,  5.49238341e+01)
( 9.00000000e+00,  5.31465394e+01)
( 1.50000000e+01,  4.35976363e+01)
( 2.80000000e+01,  3.51695047e+01)};

\nextgroupplot[yticklabels={,,}, xtick={3, 9, 15, 28}, ymin=0, xlabel={$N$ (VPD)}, ymax=100]
\addplot [black, solid, mark options={solid, thin}, mark=*, mark size=1.5, forget plot]
coordinates {
( 3.00000000e+00,  5.06139080e+01)
( 9.00000000e+00,  5.64824832e+01)
( 1.50000000e+01,  6.17508593e+01)
( 2.80000000e+01,  8.26588772e+01)};

\addplot [blue, solid, mark options={solid, thin}, mark=square*, mark size=1.5, forget plot]
coordinates {
( 3.00000000e+00,  4.79765501e+01)
( 9.00000000e+00,  5.38956431e+01)
( 1.50000000e+01,  5.24082714e+01)
( 2.80000000e+01,  6.32067097e+01)};

\addplot [red, solid, mark options={solid, thin}, mark=triangle*, mark size=1.5, forget plot]
coordinates {
( 3.00000000e+00,  5.37477253e+01)
( 9.00000000e+00,  5.77605296e+01)
( 1.50000000e+01,  5.33570702e+01)
( 2.80000000e+01,  5.49693365e+01)};

\end{groupplot}\end{tikzpicture}
 \caption{WSS reconstruction error using SBI (\textit{top row}) and MRI
 (\textit{bottom row}) as a function of MRI grid resolution for various
 noise levels (\textit{columns}) and Reynolds numbers (\textit{lines}).
 Legend: $\mathrm{Re}=100$ (\ref{line:err_sbi_vpd0}),
         $\mathrm{Re}=500$ (\ref{line:err_sbi_vpd1}),
         $\mathrm{Re}=1000$ (\ref{line:err_sbi_vpd2}).} 
 \label{fig:sten_err_wss_vpd}
\end{figure}
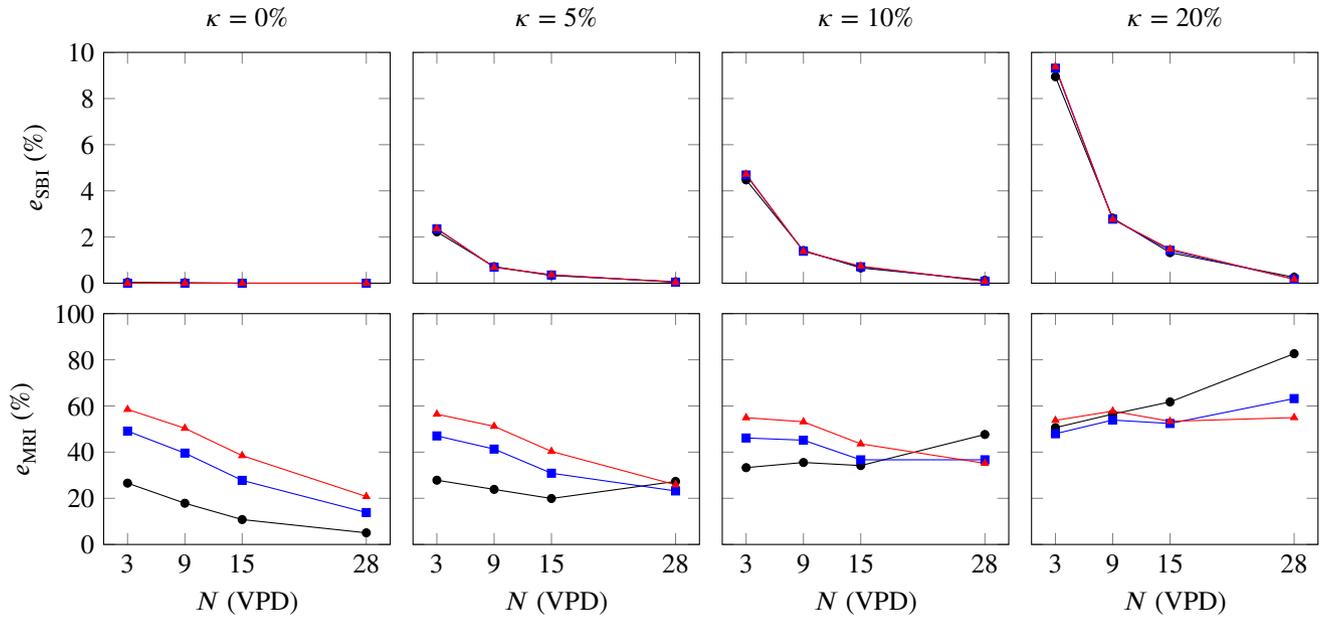

First, we focus on the relationship between the WSS reconstruction
error and resolution of the MRI voxel grid for various noise levels
and Reynolds numbers (Figure~\ref{fig:sten_err_wss_vpd}). The error
in the WSS reconstruction from MRI trends toward zero, albeit slowly,
as the voxel grid is refined, which confirms the MRI WSS approach and
implementation. Moreover, in the low noise setting,
the WSS reconstruction error decreases toward a minimum (non-zero)
value as the voxel grid is refined. However, for higher noise
levels, additional MRI voxels can degrade the WSS reconstruction
accuracy because it operates directly on the noisy MRI data and
the length scale over which the noise (fixed magnitude of $\kappa$)
varies decreases as the grid is refined, i.e., the noise
varies more rapidly in the spatial domain.
On the other hand, the WSS reconstruction from SBI is nearly
exact in the no-noise setting and the error decreases monotonically
as the voxel grid is refined for all Reynolds numbers and noise levels.
While the overall error in the WSS reconstruction using SBI does increase
with noise, additional MRI voxels do not limit or degrade the approximation
as seen in the MRI-based reconstruction. This can be attributed to the fact that
the WSS reconstruction using SBI does not \textit{directly} operate on the
noisy data; rather, the noisy data is used to reconstruct a noise-free
CFD velocity field, which is used to compute the WSS. Therefore, the
noise-sensitive operations, e.g., differentiation, are only applied
to the noise-free SBI flow field, which leads to an approximation that
is robust to noise in the MRI data. Finally, we observe that in the
$3-5$ VPD regime, the resolution commonly available for infant patients,
the WSS reconstruction using SBI is significantly more accurate than using
MRI alone. Even in the highest noise scenario ($\kappa=20\%$), the SBI
WSS reconstruction error is less than $10\%$, whereas the MRI WSS
reconstruction error is about $50\%$.

\begin{figure}
 \centering
 \begin{tikzpicture}
\begin{groupplot} [
width=0.3\textwidth,
group style={group size = 4 by 2, horizontal sep = 0.3cm, vertical sep = 0.4cm}]
\nextgroupplot[xtick={100, 500, 1000}, ymin=0, ymax=10, ylabel={$e_{\mathrm{SBI}}$ (\%)}, title={$N=3$ VPD}, xticklabels={,,}]
\addplot [black, solid, mark options={solid, thin}, mark=*, mark size=1.5]
coordinates {
( 1.00000000e+02,  3.55486846e-02)
( 5.00000000e+02,  5.72782673e-04)
( 1.00000000e+03,  8.09942723e-04)};\label{line:err_sbi_Re0}

\addplot [blue, solid, mark options={solid, thin}, mark=square*, mark size=1.5]
coordinates {
( 1.00000000e+02,  2.22485149e+00)
( 5.00000000e+02,  2.35051306e+00)
( 1.00000000e+03,  2.36878840e+00)};\label{line:err_sbi_Re1}

\addplot [red, solid, mark options={solid, thin}, mark=triangle*, mark size=1.5]
coordinates {
( 1.00000000e+02,  4.47509452e+00)
( 5.00000000e+02,  4.68653267e+00)
( 1.00000000e+03,  4.72260732e+00)};\label{line:err_sbi_Re2}

\addplot [magenta, solid, mark options={solid, thin}, mark=diamond*, mark size=1.5]
coordinates {
( 1.00000000e+02,  6.71519886e+00)
( 5.00000000e+02,  7.00734765e+00)
( 1.00000000e+03,  7.06046404e+00)};\label{line:err_sbi_Re3}

\addplot [green, solid, mark options={solid, thin}, mark=pentagon*, mark size=1.5]
coordinates {
( 1.00000000e+02,  8.94502468e+00)
( 5.00000000e+02,  9.31281610e+00)
( 1.00000000e+03,  9.38217183e+00)};\label{line:err_sbi_Re4}

\nextgroupplot[yticklabels={,,}, xtick={100, 500, 1000}, ymax=10, title={$N=9$ VPD}, ymin=0, xticklabels={,,}]
\addplot [black, solid, mark options={solid, thin}, mark=*, mark size=1.5, forget plot]
coordinates {
( 1.00000000e+02,  2.11432050e-02)
( 5.00000000e+02,  1.41276278e-05)
( 1.00000000e+03,  1.37969706e-05)};

\addplot [blue, solid, mark options={solid, thin}, mark=square*, mark size=1.5, forget plot]
coordinates {
( 1.00000000e+02,  7.19548739e-01)
( 5.00000000e+02,  6.91964759e-01)
( 1.00000000e+03,  6.90204685e-01)};

\addplot [red, solid, mark options={solid, thin}, mark=triangle*, mark size=1.5, forget plot]
coordinates {
( 1.00000000e+02,  1.41869455e+00)
( 5.00000000e+02,  1.38500199e+00)
( 1.00000000e+03,  1.38190010e+00)};

\addplot [magenta, solid, mark options={solid, thin}, mark=diamond*, mark size=1.5, forget plot]
coordinates {
( 1.00000000e+02,  2.11910222e+00)
( 5.00000000e+02,  2.07912483e+00)
( 1.00000000e+03,  2.07457007e+00)};

\addplot [green, solid, mark options={solid, thin}, mark=pentagon*, mark size=1.5, forget plot]
coordinates {
( 1.00000000e+02,  2.81972424e+00)
( 5.00000000e+02,  2.77433231e+00)
( 1.00000000e+03,  2.76818846e+00)};

\nextgroupplot[yticklabels={,,}, xtick={100, 500, 1000}, ymax=10, title={$N=15$ VPD}, ymin=0, xticklabels={,,}]
\addplot [black, solid, mark options={solid, thin}, mark=*, mark size=1.5, forget plot]
coordinates {
( 1.00000000e+02,  3.96954920e-05)
( 5.00000000e+02,  1.64484550e-05)
( 1.00000000e+03,  2.06027349e-05)};

\addplot [blue, solid, mark options={solid, thin}, mark=square*, mark size=1.5, forget plot]
coordinates {
( 1.00000000e+02,  3.30628594e-01)
( 5.00000000e+02,  3.55280957e-01)
( 1.00000000e+03,  3.69171269e-01)};

\addplot [red, solid, mark options={solid, thin}, mark=triangle*, mark size=1.5, forget plot]
coordinates {
( 1.00000000e+02,  6.61457667e-01)
( 5.00000000e+02,  7.11002174e-01)
( 1.00000000e+03,  7.38778108e-01)};

\addplot [magenta, solid, mark options={solid, thin}, mark=diamond*, mark size=1.5, forget plot]
coordinates {
( 1.00000000e+02,  9.92448012e-01)
( 5.00000000e+02,  1.06714676e+00)
( 1.00000000e+03,  1.10884105e+00)};

\addplot [green, solid, mark options={solid, thin}, mark=pentagon*, mark size=1.5, forget plot]
coordinates {
( 1.00000000e+02,  1.32360011e+00)
( 5.00000000e+02,  1.42371425e+00)
( 1.00000000e+03,  1.47935571e+00)};

\nextgroupplot[yticklabels={,,}, xtick={100, 500, 1000}, ymax=10, title={$N=28$ VPD}, ymin=0, xticklabels={,,}]
\addplot [black, solid, mark options={solid, thin}, mark=*, mark size=1.5, forget plot]
coordinates {
( 1.00000000e+02,  1.92496800e-04)
( 5.00000000e+02,  2.34543308e-05)
( 1.00000000e+03,  5.83797456e-06)};

\addplot [blue, solid, mark options={solid, thin}, mark=square*, mark size=1.5, forget plot]
coordinates {
( 1.00000000e+02,  6.46523112e-02)
( 5.00000000e+02,  4.41273842e-02)
( 1.00000000e+03,  4.26644615e-02)};

\addplot [red, solid, mark options={solid, thin}, mark=triangle*, mark size=1.5, forget plot]
coordinates {
( 1.00000000e+02,  1.29114460e-01)
( 5.00000000e+02,  8.82311725e-02)
( 1.00000000e+03,  8.53286940e-02)};

\addplot [magenta, solid, mark options={solid, thin}, mark=diamond*, mark size=1.5, forget plot]
coordinates {
( 1.00000000e+02,  1.93578926e-01)
( 5.00000000e+02,  1.32334825e-01)
( 1.00000000e+03,  1.27986863e-01)};

\addplot [green, solid, mark options={solid, thin}, mark=pentagon*, mark size=1.5, forget plot]
coordinates {
( 1.00000000e+02,  2.58045694e-01)
( 5.00000000e+02,  1.76438350e-01)
( 1.00000000e+03,  1.70638973e-01)};

\nextgroupplot[xtick={100, 500, 1000}, ymin=0, xlabel={$\mathrm{Re}$}, ymax=100, ylabel={$e_{\mathrm{MRI}}$ (\%)}]
\addplot [black, solid, mark options={solid, thin}, mark=*, mark size=1.5]
coordinates {
( 1.00000000e+02,  2.65715241e+01)
( 5.00000000e+02,  4.91077861e+01)
( 1.00000000e+03,  5.85280535e+01)};\label{line:err_mri_Re0}

\addplot [blue, solid, mark options={solid, thin}, mark=square*, mark size=1.5]
coordinates {
( 1.00000000e+02,  2.78437783e+01)
( 5.00000000e+02,  4.70431007e+01)
( 1.00000000e+03,  5.64441610e+01)};\label{line:err_mri_Re1}

\addplot [red, solid, mark options={solid, thin}, mark=triangle*, mark size=1.5]
coordinates {
( 1.00000000e+02,  3.32812207e+01)
( 5.00000000e+02,  4.61364859e+01)
( 1.00000000e+03,  5.49238341e+01)};\label{line:err_mri_Re2}

\addplot [magenta, solid, mark options={solid, thin}, mark=diamond*, mark size=1.5]
coordinates {
( 1.00000000e+02,  4.12690790e+01)
( 5.00000000e+02,  4.64557927e+01)
( 1.00000000e+03,  5.40146811e+01)};\label{line:err_mri_Re3}

\addplot [green, solid, mark options={solid, thin}, mark=pentagon*, mark size=1.5]
coordinates {
( 1.00000000e+02,  5.06139080e+01)
( 5.00000000e+02,  4.79765501e+01)
( 1.00000000e+03,  5.37477253e+01)};\label{line:err_mri_Re4}

\nextgroupplot[yticklabels={,,}, xtick={100, 500, 1000}, ymin=0, xlabel={$\mathrm{Re}$}, ymax=100]
\addplot [black, solid, mark options={solid, thin}, mark=*, mark size=1.5, forget plot]
coordinates {
( 1.00000000e+02,  1.78628927e+01)
( 5.00000000e+02,  3.95811298e+01)
( 1.00000000e+03,  5.03418769e+01)};

\addplot [blue, solid, mark options={solid, thin}, mark=square*, mark size=1.5, forget plot]
coordinates {
( 1.00000000e+02,  2.38641502e+01)
( 5.00000000e+02,  4.13301593e+01)
( 1.00000000e+03,  5.11903870e+01)};

\addplot [red, solid, mark options={solid, thin}, mark=triangle*, mark size=1.5, forget plot]
coordinates {
( 1.00000000e+02,  3.55183830e+01)
( 5.00000000e+02,  4.51786292e+01)
( 1.00000000e+03,  5.31465394e+01)};

\addplot [magenta, solid, mark options={solid, thin}, mark=diamond*, mark size=1.5, forget plot]
coordinates {
( 1.00000000e+02,  4.52149344e+01)
( 5.00000000e+02,  4.95189670e+01)
( 1.00000000e+03,  5.56024814e+01)};

\addplot [green, solid, mark options={solid, thin}, mark=pentagon*, mark size=1.5, forget plot]
coordinates {
( 1.00000000e+02,  5.64824832e+01)
( 5.00000000e+02,  5.38956431e+01)
( 1.00000000e+03,  5.77605296e+01)};

\nextgroupplot[yticklabels={,,}, xtick={100, 500, 1000}, ymin=0, xlabel={$\mathrm{Re}$}, ymax=100]
\addplot [black, solid, mark options={solid, thin}, mark=*, mark size=1.5, forget plot]
coordinates {
( 1.00000000e+02,  1.07892147e+01)
( 5.00000000e+02,  2.77682548e+01)
( 1.00000000e+03,  3.84932763e+01)};

\addplot [blue, solid, mark options={solid, thin}, mark=square*, mark size=1.5, forget plot]
coordinates {
( 1.00000000e+02,  1.99183403e+01)
( 5.00000000e+02,  3.08691039e+01)
( 1.00000000e+03,  4.03305203e+01)};

\addplot [red, solid, mark options={solid, thin}, mark=triangle*, mark size=1.5, forget plot]
coordinates {
( 1.00000000e+02,  3.41795120e+01)
( 5.00000000e+02,  3.66709606e+01)
( 1.00000000e+03,  4.35976363e+01)};

\addplot [magenta, solid, mark options={solid, thin}, mark=diamond*, mark size=1.5, forget plot]
coordinates {
( 1.00000000e+02,  4.85691600e+01)
( 5.00000000e+02,  4.41166207e+01)
( 1.00000000e+03,  4.80430180e+01)};

\addplot [green, solid, mark options={solid, thin}, mark=pentagon*, mark size=1.5, forget plot]
coordinates {
( 1.00000000e+02,  6.17508593e+01)
( 5.00000000e+02,  5.24082714e+01)
( 1.00000000e+03,  5.33570702e+01)};

\nextgroupplot[yticklabels={,,}, xtick={100, 500, 1000}, ymin=0, xlabel={$\mathrm{Re}$}, ymax=100]
\addplot [black, solid, mark options={solid, thin}, mark=*, mark size=1.5, forget plot]
coordinates {
( 1.00000000e+02,  5.05153628e+00)
( 5.00000000e+02,  1.38287456e+01)
( 1.00000000e+03,  2.07782254e+01)};

\addplot [blue, solid, mark options={solid, thin}, mark=square*, mark size=1.5, forget plot]
coordinates {
( 1.00000000e+02,  2.73248578e+01)
( 5.00000000e+02,  2.31931672e+01)
( 1.00000000e+03,  2.58547775e+01)};

\addplot [red, solid, mark options={solid, thin}, mark=triangle*, mark size=1.5, forget plot]
coordinates {
( 1.00000000e+02,  4.76494314e+01)
( 5.00000000e+02,  3.66924359e+01)
( 1.00000000e+03,  3.51695047e+01)};

\addplot [magenta, solid, mark options={solid, thin}, mark=diamond*, mark size=1.5, forget plot]
coordinates {
( 1.00000000e+02,  6.59568901e+01)
( 5.00000000e+02,  4.95818300e+01)
( 1.00000000e+03,  4.47011917e+01)};

\addplot [green, solid, mark options={solid, thin}, mark=pentagon*, mark size=1.5, forget plot]
coordinates {
( 1.00000000e+02,  8.26588772e+01)
( 5.00000000e+02,  6.32067097e+01)
( 1.00000000e+03,  5.49693365e+01)};

\end{groupplot}\end{tikzpicture}
 \caption{WSS reconstruction error using SBI (\textit{top row}) and MRI
 (\textit{bottom row}) as a function of Reynolds number for various
 MRI grid resolutions (\textit{columns}) and noise levels (\textit{lines}).
 Legend: $\kappa=0\%$ (\ref{line:err_sbi_Re0}),
         $\kappa=5\%$ (\ref{line:err_sbi_Re1}),
         $\kappa=10\%$ (\ref{line:err_sbi_Re2}),
         $\kappa=15\%$ (\ref{line:err_sbi_Re3}),
         $\kappa=20\%$ (\ref{line:err_sbi_Re4}).}
 \label{fig:sten_err_wss_Re}
\end{figure}
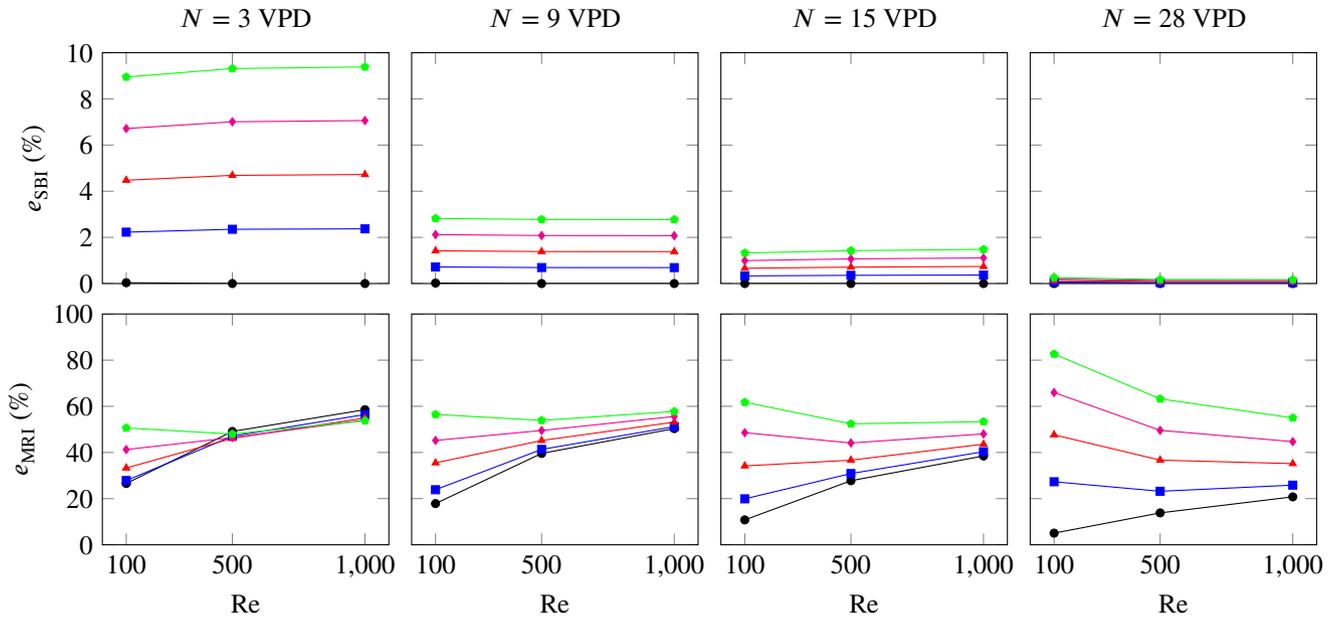

Next, we investigate the relationship between the WSS reconstruction
error and Reynolds number for various noise levels and MRI voxel grids
(Figure~\ref{fig:sten_err_wss_Re}). The accuracy of the WSS reconstruction
directly from the MRI degrades as the Reynolds number increases, which
agrees with other studies \cite{quant_wss_numerical}. The main exceptions
come from configurations
with high noise and many voxels per diameter where the error is already quite
large (above $40\%$); in these cases, the WSS reconstruction error can decrease
somewhat with increasing Reynolds number. On the other hand, WSS reconstruction
from SBI is insensitive to Reynolds number, 
i.e., the range in the WSS reconstruction error is
less than $1\%$ from $\mathrm{Re}=100$ to $\mathrm{Re}=1000$ for all noise
levels and decreases as the MRI grid is refined.

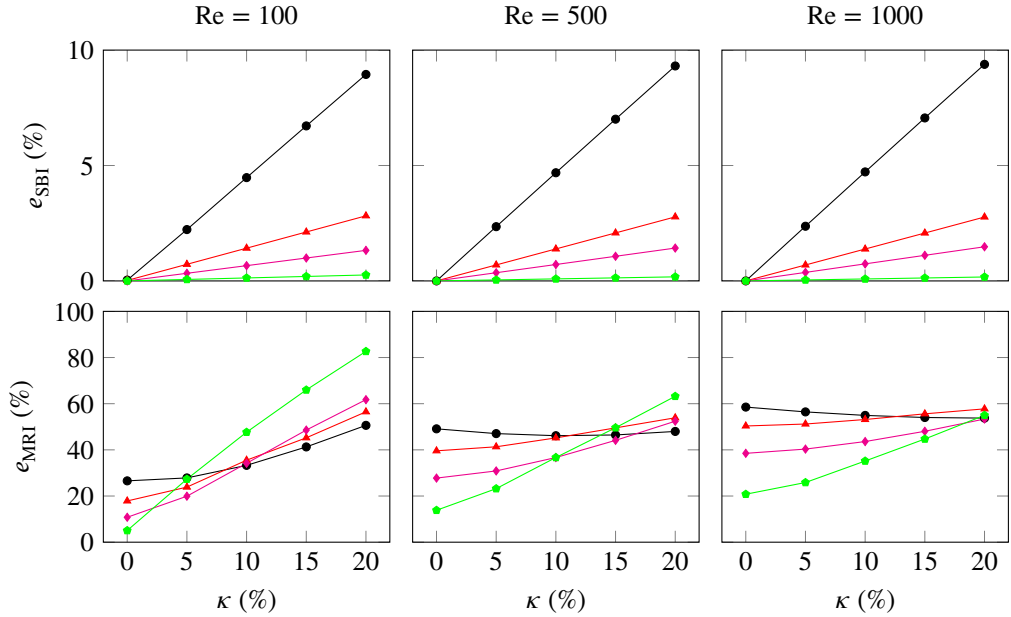
\begin{figure}
 \centering
 \begin{tikzpicture}
\begin{groupplot} [
width=0.3\textwidth,
group style={group size = 3 by 2, horizontal sep = 0.3cm, vertical sep = 0.4cm}]
\nextgroupplot[xtick={0, 5, 10, 15, 20}, ytick={0,5,10}, ymax=10, title={$\mathrm{Re}=100$}, ymin=0, ylabel={$e_{\mathrm{SBI}}$ (\%)}, xticklabels={,,}]
\addplot [black, solid, mark options={solid, thin}, mark=*, mark size=1.5]
coordinates {
( 0.00000000e+00,  3.55486846e-02)
( 5.00000000e+00,  2.22485149e+00)
( 1.00000000e+01,  4.47509452e+00)
( 1.50000000e+01,  6.71519886e+00)
( 2.00000000e+01,  8.94502468e+00)};\label{line:err_sbi_noise0}

\addplot [red, solid, mark options={solid, thin}, mark=triangle*, mark size=1.5]
coordinates {
( 0.00000000e+00,  2.11432050e-02)
( 5.00000000e+00,  7.19548739e-01)
( 1.00000000e+01,  1.41869455e+00)
( 1.50000000e+01,  2.11910222e+00)
( 2.00000000e+01,  2.81972424e+00)};\label{line:err_sbi_noise2}

\addplot [magenta, solid, mark options={solid, thin}, mark=diamond*, mark size=1.5]
coordinates {
( 0.00000000e+00,  3.96954920e-05)
( 5.00000000e+00,  3.30628594e-01)
( 1.00000000e+01,  6.61457667e-01)
( 1.50000000e+01,  9.92448012e-01)
( 2.00000000e+01,  1.32360011e+00)};\label{line:err_sbi_noise3}

\addplot [green, solid, mark options={solid, thin}, mark=pentagon*, mark size=1.5]
coordinates {
( 0.00000000e+00,  1.92496800e-04)
( 5.00000000e+00,  6.46523112e-02)
( 1.00000000e+01,  1.29114460e-01)
( 1.50000000e+01,  1.93578926e-01)
( 2.00000000e+01,  2.58045694e-01)};\label{line:err_sbi_noise4}

\nextgroupplot[yticklabels={,,}, xtick={0, 5, 10, 15, 20}, ytick={0,5,10}, ymax=10, title={$\mathrm{Re}=500$}, ymin=0, xticklabels={,,}]
\addplot [black, solid, mark options={solid, thin}, mark=*, mark size=1.5, forget plot]
coordinates {
( 0.00000000e+00,  5.72782673e-04)
( 5.00000000e+00,  2.35051306e+00)
( 1.00000000e+01,  4.68653267e+00)
( 1.50000000e+01,  7.00734765e+00)
( 2.00000000e+01,  9.31281610e+00)};

\addplot [red, solid, mark options={solid, thin}, mark=triangle*, mark size=1.5, forget plot]
coordinates {
( 0.00000000e+00,  1.41276278e-05)
( 5.00000000e+00,  6.91964759e-01)
( 1.00000000e+01,  1.38500199e+00)
( 1.50000000e+01,  2.07912483e+00)
( 2.00000000e+01,  2.77433231e+00)};

\addplot [magenta, solid, mark options={solid, thin}, mark=diamond*, mark size=1.5, forget plot]
coordinates {
( 0.00000000e+00,  1.64484550e-05)
( 5.00000000e+00,  3.55280957e-01)
( 1.00000000e+01,  7.11002174e-01)
( 1.50000000e+01,  1.06714676e+00)
( 2.00000000e+01,  1.42371425e+00)};

\addplot [green, solid, mark options={solid, thin}, mark=pentagon*, mark size=1.5, forget plot]
coordinates {
( 0.00000000e+00,  2.34543308e-05)
( 5.00000000e+00,  4.41273842e-02)
( 1.00000000e+01,  8.82311725e-02)
( 1.50000000e+01,  1.32334825e-01)
( 2.00000000e+01,  1.76438350e-01)};

\nextgroupplot[yticklabels={,,}, xtick={0, 5, 10, 15, 20}, ytick={0,5,10}, ymax=10, title={$\mathrm{Re}=1000$}, ymin=0, xticklabels={,,}]
\addplot [black, solid, mark options={solid, thin}, mark=*, mark size=1.5, forget plot]
coordinates {
( 0.00000000e+00,  8.09942723e-04)
( 5.00000000e+00,  2.36878840e+00)
( 1.00000000e+01,  4.72260732e+00)
( 1.50000000e+01,  7.06046404e+00)
( 2.00000000e+01,  9.38217183e+00)};

\addplot [red, solid, mark options={solid, thin}, mark=triangle*, mark size=1.5, forget plot]
coordinates {
( 0.00000000e+00,  1.37969706e-05)
( 5.00000000e+00,  6.90204685e-01)
( 1.00000000e+01,  1.38190010e+00)
( 1.50000000e+01,  2.07457007e+00)
( 2.00000000e+01,  2.76818846e+00)};

\addplot [magenta, solid, mark options={solid, thin}, mark=diamond*, mark size=1.5, forget plot]
coordinates {
( 0.00000000e+00,  2.06027349e-05)
( 5.00000000e+00,  3.69171269e-01)
( 1.00000000e+01,  7.38778108e-01)
( 1.50000000e+01,  1.10884105e+00)
( 2.00000000e+01,  1.47935571e+00)};

\addplot [green, solid, mark options={solid, thin}, mark=pentagon*, mark size=1.5, forget plot]
coordinates {
( 0.00000000e+00,  5.83797456e-06)
( 5.00000000e+00,  4.26644615e-02)
( 1.00000000e+01,  8.53286940e-02)
( 1.50000000e+01,  1.27986863e-01)
( 2.00000000e+01,  1.70638973e-01)};

\nextgroupplot[xtick={0, 5, 10, 15, 20}, ymin=0, xlabel={$\kappa$ (\%)}, ymax=100, ylabel={$e_{\mathrm{MRI}}$ (\%)}]
\addplot [black, solid, mark options={solid, thin}, mark=*, mark size=1.5, forget plot]
coordinates {
( 0.00000000e+00,  2.65715241e+01)
( 5.00000000e+00,  2.78437783e+01)
( 1.00000000e+01,  3.32812207e+01)
( 1.50000000e+01,  4.12690790e+01)
( 2.00000000e+01,  5.06139080e+01)};

\addplot [red, solid, mark options={solid, thin}, mark=triangle*, mark size=1.5, forget plot]
coordinates {
( 0.00000000e+00,  1.78628927e+01)
( 5.00000000e+00,  2.38641502e+01)
( 1.00000000e+01,  3.55183830e+01)
( 1.50000000e+01,  4.52149344e+01)
( 2.00000000e+01,  5.64824832e+01)};

\addplot [magenta, solid, mark options={solid, thin}, mark=diamond*, mark size=1.5, forget plot]
coordinates {
( 0.00000000e+00,  1.07892147e+01)
( 5.00000000e+00,  1.99183403e+01)
( 1.00000000e+01,  3.41795120e+01)
( 1.50000000e+01,  4.85691600e+01)
( 2.00000000e+01,  6.17508593e+01)};

\addplot [green, solid, mark options={solid, thin}, mark=pentagon*, mark size=1.5, forget plot]
coordinates {
( 0.00000000e+00,  5.05153628e+00)
( 5.00000000e+00,  2.73248578e+01)
( 1.00000000e+01,  4.76494314e+01)
( 1.50000000e+01,  6.59568901e+01)
( 2.00000000e+01,  8.26588772e+01)};

\nextgroupplot[yticklabels={,,}, xtick={0, 5, 10, 15, 20}, ymin=0, xlabel={$\kappa$ (\%)}, ymax=100]
\addplot [black, solid, mark options={solid, thin}, mark=*, mark size=1.5, forget plot]
coordinates {
( 0.00000000e+00,  4.91077861e+01)
( 5.00000000e+00,  4.70431007e+01)
( 1.00000000e+01,  4.61364859e+01)
( 1.50000000e+01,  4.64557927e+01)
( 2.00000000e+01,  4.79765501e+01)};

\addplot [red, solid, mark options={solid, thin}, mark=triangle*, mark size=1.5, forget plot]
coordinates {
( 0.00000000e+00,  3.95811298e+01)
( 5.00000000e+00,  4.13301593e+01)
( 1.00000000e+01,  4.51786292e+01)
( 1.50000000e+01,  4.95189670e+01)
( 2.00000000e+01,  5.38956431e+01)};

\addplot [magenta, solid, mark options={solid, thin}, mark=diamond*, mark size=1.5, forget plot]
coordinates {
( 0.00000000e+00,  2.77682548e+01)
( 5.00000000e+00,  3.08691039e+01)
( 1.00000000e+01,  3.66709606e+01)
( 1.50000000e+01,  4.41166207e+01)
( 2.00000000e+01,  5.24082714e+01)};

\addplot [green, solid, mark options={solid, thin}, mark=pentagon*, mark size=1.5, forget plot]
coordinates {
( 0.00000000e+00,  1.38287456e+01)
( 5.00000000e+00,  2.31931672e+01)
( 1.00000000e+01,  3.66924359e+01)
( 1.50000000e+01,  4.95818300e+01)
( 2.00000000e+01,  6.32067097e+01)};

\nextgroupplot[yticklabels={,,}, xtick={0, 5, 10, 15, 20}, ymin=0, xlabel={$\kappa$ (\%)}, ymax=100]
\addplot [black, solid, mark options={solid, thin}, mark=*, mark size=1.5, forget plot]
coordinates {
( 0.00000000e+00,  5.85280535e+01)
( 5.00000000e+00,  5.64441610e+01)
( 1.00000000e+01,  5.49238341e+01)
( 1.50000000e+01,  5.40146811e+01)
( 2.00000000e+01,  5.37477253e+01)};

\addplot [red, solid, mark options={solid, thin}, mark=triangle*, mark size=1.5, forget plot]
coordinates {
( 0.00000000e+00,  5.03418769e+01)
( 5.00000000e+00,  5.11903870e+01)
( 1.00000000e+01,  5.31465394e+01)
( 1.50000000e+01,  5.56024814e+01)
( 2.00000000e+01,  5.77605296e+01)};

\addplot [magenta, solid, mark options={solid, thin}, mark=diamond*, mark size=1.5, forget plot]
coordinates {
( 0.00000000e+00,  3.84932763e+01)
( 5.00000000e+00,  4.03305203e+01)
( 1.00000000e+01,  4.35976363e+01)
( 1.50000000e+01,  4.80430180e+01)
( 2.00000000e+01,  5.33570702e+01)};

\addplot [green, solid, mark options={solid, thin}, mark=pentagon*, mark size=1.5, forget plot]
coordinates {
( 0.00000000e+00,  2.07782254e+01)
( 5.00000000e+00,  2.58547775e+01)
( 1.00000000e+01,  3.51695047e+01)
( 1.50000000e+01,  4.47011917e+01)
( 2.00000000e+01,  5.49693365e+01)};

\end{groupplot}\end{tikzpicture}
 \caption{WSS reconstruction error using SBI (\textit{top row}) and MRI
 (\textit{bottom row}) as a function of noise level for various
 Reynolds numbers (\textit{columns}) and MRI grid resolutions (\textit{lines}).
 Legend: $N=3$ VPD (\ref{line:err_sbi_noise0}),
         $N=9$ VPD (\ref{line:err_sbi_noise2}),
         $N=15$ VPD (\ref{line:err_sbi_noise3}),
         $N=28$ VPD (\ref{line:err_sbi_noise4}).}
    \label{fig:sten_err_wss_noise}
\end{figure}

Next, we investigate the relationship between the WSS reconstruction
error and noise level for various Reynolds numbers and MRI voxel grids
(Figure~\ref{fig:sten_err_wss_noise}). For the MRI-based WSS reconstruction,
the error increases with noise except for the coarsest voxel grid
($N=3$ VPD) where the error can slightly decrease as noise increases for higher
Reynolds numbers (situations where the error is already quite large, at least
$50\%$). As the MRI grid is refined, the rate at which the WSS error increases
with respect to noise accelerates due to the decreasing length scale over which the noise varies.
The error in the WSS reconstruction using SBI increases linearly with the
noise level with a decreasing slope as the MRI resolution increases (opposite
of the trend observed with MRI-only WSS reconstruction).
These observations hold across all Reynolds numbers considered
$\mathrm{Re}\in\{100,500,1000\}$.

\begin{figure}
 \centering
 \input{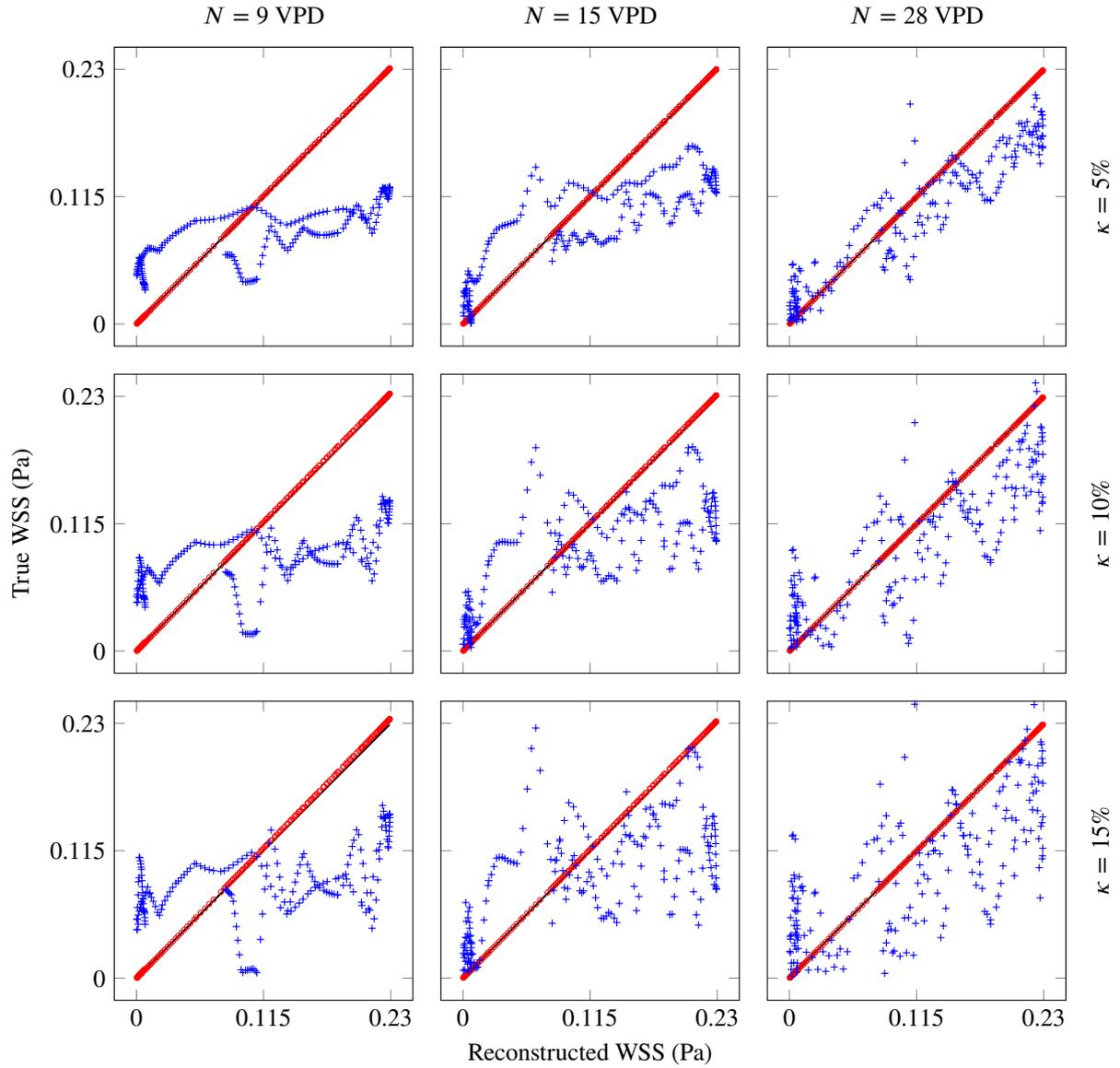}
 \caption{Scatter plot of the true WSS vs. the reconstructed WSS using
 SBI (\ref{line:wss_dist_sbi00}) and MRI (\ref{line:wss_dist_mri00}) for
 $200$ points along the upper wall intersected with the MRI domain ($\Gamma$).
 The \textit{rows} correspond to noise levels $\kappa=5\%,10\%,15\%$ and
 the \textit{columns} correspond to MRI resolutions $N=9,15,28$ VPD.
 Tight clustering around the line of identity (\ref{line:wss_dist_exact00})
 indicates accurate WSS reconstruction along $\Gamma$.}
 \label{fig:sten_wss_dist}
\end{figure}

To this point, we have aggregated the entire WSS distribution into a scalar
to compare MRI- and SBI-based WSS reconstruction across numerous scenarios. 
We take a closer look at the entire WSS distribution as a scatter plot,
whereby the actual value of WSS is plotted against the reconstructed WSS for
a number of points along the wall $x\in\Gamma$ (Figure~\ref{fig:sten_wss_dist});
tight clustering around the line of identity implies the reconstruction is
accurate and reliable. Because of the relatively weak dependence on Reynolds
number (Figure~\ref{fig:sten_err_wss_Re}), we fix $\mathrm{Re} = 1000$ and
vary the noise
$\kappa\in\{5\%,10\%,15\%\}$ and MRI resolution $N\in\{9,15, 28\}$ VPD.
The SBI-based WSS reconstruction lies tightly clustered near the line
of identity, whereas the MRI-based WSS reconstruction varies significantly
from the line, particularly in the high-noise configurations.

\subsection{Verification with aorta geometry}
A final set of numerical experiments was conducted using the coarctated
aorta to demonstrate the key findings from the previous sections generalize
to the more complex test case. To limit the parameter space to explore,
we do not study the impact of mesh resolution on the SBI reconstruction.
Furthermore, given the relative insensitivity of SBI to Reynolds number,
we only consider a limited sampling of Reynolds numbers, i.e.,
$\mathrm{Re}\in\{100, 1000\}$. We vary the noise $\kappa\in\{0, 5, 10, 20\}$
and MRI voxel grid $N\in\{3, 5, 10\}$ VPD; the finest MRI grid is restricted
to $N=10$ VPD (corresponds to $60\times 90$ voxel grid) because the MRI domain
covers a larger region compared to the stenotic vessel.

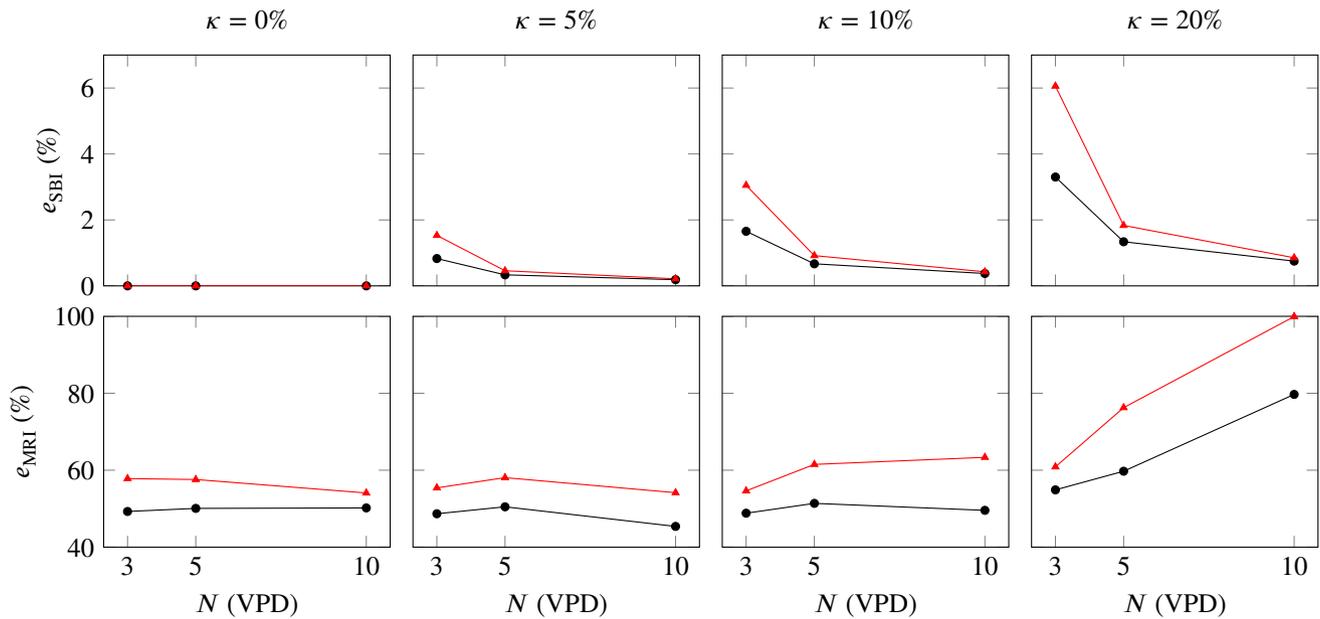
\begin{figure}
 \centering
 \begin{tikzpicture}
\begin{groupplot} [
width=0.3\textwidth,
group style={group size = 4 by 2, horizontal sep = 0.3cm, vertical sep = 0.4cm}]
\nextgroupplot[xtick={3, 5, 10}, ymin=0, ymax=7, ylabel={$e_{\mathrm{SBI}}$ (\%)}, title={$\kappa=0\%$}, xticklabels={,,}]
\addplot [black, solid, mark options={solid, thin}, mark=*, mark size=1.5]
coordinates {
( 3.00000000e+00,  3.35266163e-04)
( 5.00000000e+00,  1.33053288e-04)
( 1.00000000e+01,  3.34920739e-05)};\label{line:err_aorta_sbi_vpd0}

\addplot [red, solid, mark options={solid, thin}, mark=triangle*, mark size=1.5]
coordinates {
( 3.00000000e+00,  3.10178650e-06)
( 5.00000000e+00,  1.27440281e-06)
( 1.00000000e+01,  3.07814548e-07)};\label{line:err_aorta_sbi_vpd2}

\nextgroupplot[yticklabels={,,}, xtick={3, 5, 10}, ymax=7, title={$\kappa=5\%$}, ymin=0, xticklabels={,,}]
\addplot [black, solid, mark options={solid, thin}, mark=*, mark size=1.5, forget plot]
coordinates {
( 3.00000000e+00,  8.28785044e-01)
( 5.00000000e+00,  3.34642952e-01)
( 1.00000000e+01,  1.87697228e-01)};

\addplot [red, solid, mark options={solid, thin}, mark=triangle*, mark size=1.5, forget plot]
coordinates {
( 3.00000000e+00,  1.53091772e+00)
( 5.00000000e+00,  4.59039191e-01)
( 1.00000000e+01,  2.14396861e-01)};

\nextgroupplot[yticklabels={,,}, xtick={3, 5, 10}, ymax=7, title={$\kappa=10\%$}, ymin=0, xticklabels={,,}]
\addplot [black, solid, mark options={solid, thin}, mark=*, mark size=1.5, forget plot]
coordinates {
( 3.00000000e+00,  1.65495715e+00)
( 5.00000000e+00,  6.69272793e-01)
( 1.00000000e+01,  3.75426635e-01)};

\addplot [red, solid, mark options={solid, thin}, mark=triangle*, mark size=1.5, forget plot]
coordinates {
( 3.00000000e+00,  3.05017086e+00)
( 5.00000000e+00,  9.17523856e-01)
( 1.00000000e+01,  4.27846301e-01)};

\nextgroupplot[yticklabels={,,}, xtick={3, 5, 10}, ymax=7, title={$\kappa=20\%$}, ymin=0, xticklabels={,,}]
\addplot [black, solid, mark options={solid, thin}, mark=*, mark size=1.5, forget plot]
coordinates {
( 3.00000000e+00,  3.30042217e+00)
( 5.00000000e+00,  1.33809588e+00)
( 1.00000000e+01,  7.50959854e-01)};

\addplot [red, solid, mark options={solid, thin}, mark=triangle*, mark size=1.5, forget plot]
coordinates {
( 3.00000000e+00,  6.05798417e+00)
( 5.00000000e+00,  1.83344958e+00)
( 1.00000000e+01,  8.52155239e-01)};

\nextgroupplot[xtick={3, 5, 10}, ymin=40, xlabel={$N$ (VPD)}, ymax=100, ylabel={$e_{\mathrm{MRI}}$ (\%)}]
\addplot [black, solid, mark options={solid, thin}, mark=*, mark size=1.5, forget plot]
coordinates {
( 3.00000000e+00,  4.92727793e+01)
( 5.00000000e+00,  5.00731302e+01)
( 1.00000000e+01,  5.01915016e+01)};

\addplot [red, solid, mark options={solid, thin}, mark=triangle*, mark size=1.5, forget plot]
coordinates {
( 3.00000000e+00,  5.78198694e+01)
( 5.00000000e+00,  5.75788626e+01)
( 1.00000000e+01,  5.40835180e+01)};

\nextgroupplot[yticklabels={,,}, xtick={3, 5, 10}, ymin=40, xlabel={$N$ (VPD)}, ymax=100]
\addplot [black, solid, mark options={solid, thin}, mark=*, mark size=1.5, forget plot]
coordinates {
( 3.00000000e+00,  4.86778157e+01)
( 5.00000000e+00,  5.04804451e+01)
( 1.00000000e+01,  4.53927526e+01)};

\addplot [red, solid, mark options={solid, thin}, mark=triangle*, mark size=1.5, forget plot]
coordinates {
( 3.00000000e+00,  5.54113014e+01)
( 5.00000000e+00,  5.80649285e+01)
( 1.00000000e+01,  5.41498969e+01)};

\nextgroupplot[yticklabels={,,}, xtick={3, 5, 10}, ymin=40, xlabel={$N$ (VPD)}, ymax=100]
\addplot [black, solid, mark options={solid, thin}, mark=*, mark size=1.5, forget plot]
coordinates {
( 3.00000000e+00,  4.88207029e+01)
( 5.00000000e+00,  5.13778601e+01)
( 1.00000000e+01,  4.95657434e+01)};

\addplot [red, solid, mark options={solid, thin}, mark=triangle*, mark size=1.5, forget plot]
coordinates {
( 3.00000000e+00,  5.46254465e+01)
( 5.00000000e+00,  6.15050846e+01)
( 1.00000000e+01,  6.33429779e+01)};

\nextgroupplot[yticklabels={,,}, xtick={3, 5, 10}, ymin=40, xlabel={$N$ (VPD)}, ymax=100]
\addplot [black, solid, mark options={solid, thin}, mark=*, mark size=1.5, forget plot]
coordinates {
( 3.00000000e+00,  5.48890428e+01)
( 5.00000000e+00,  5.97287299e+01)
( 1.00000000e+01,  7.96994298e+01)};

\addplot [red, solid, mark options={solid, thin}, mark=triangle*, mark size=1.5, forget plot]
coordinates {
( 3.00000000e+00,  6.08726002e+01)
( 5.00000000e+00,  7.62583799e+01)
( 1.00000000e+01,  9.99386324e+01)};

\end{groupplot}\end{tikzpicture}
 \caption{WSS reconstruction error using SBI (\textit{top row}) and MRI
 (\textit{bottom row}) as a function of MRI grid resolution for various
 noise levels (\textit{columns}) and Reynolds numbers (\textit{lines})
 for coarctated aorta test case.
 Legend: $\mathrm{Re}=100$ (\ref{line:err_aorta_sbi_vpd0}),
         $\mathrm{Re}=1000$ (\ref{line:err_aorta_sbi_vpd2}).}
 \label{fig:aorta_err_wss_vpd}
\end{figure}

The relationship between the WSS reconstruction error and resolution of the
MRI voxel grid for various noise levels and Reynolds numbers
(Figure~\ref{fig:aorta_err_wss_vpd})
is consistent with findings for the stenotic vessel. That is, the accuracy
of the WSS reconstruction from MRI improves as the MRI voxel grid is refined
for the low-noise settings, whereas it degrades in the higher noise settings.
On the other hand, the WSS reconstruction error from SBI decreases as the MRI
voxel grid is refined for all noise levels considered.

\begin{figure}
 \centering
 \begin{tikzpicture}
\begin{groupplot} [
width=0.3\textwidth,
group style={group size = 2 by 2, horizontal sep = 0.3cm, vertical sep = 0.4cm}]
\nextgroupplot[xtick={0, 5, 10, 20}, ytick={0,3,6,9,12}, ymax=7, title={$\mathrm{Re}=100$}, ymin=0, ylabel={$e_{\mathrm{SBI}}$ (\%)}, xticklabels={,,}]
\addplot [black, solid, mark options={solid, thin}, mark=*, mark size=1.5]
coordinates {
( 0.00000000e+00,  3.35266163e-04)
( 5.00000000e+00,  8.28785044e-01)
( 1.00000000e+01,  1.65495715e+00)
( 2.00000000e+01,  3.30042217e+00)};\label{line:err_aorta_sbi_noise0}

\addplot [blue, solid, mark options={solid, thin}, mark=square*, mark size=1.5]
coordinates {
( 0.00000000e+00,  1.33053288e-04)
( 5.00000000e+00,  3.34642952e-01)
( 1.00000000e+01,  6.69272793e-01)
( 2.00000000e+01,  1.33809588e+00)};\label{line:err_aorta_sbi_noise1}

\addplot [red, solid, mark options={solid, thin}, mark=triangle*, mark size=1.5]
coordinates {
( 0.00000000e+00,  3.34920739e-05)
( 5.00000000e+00,  1.87697228e-01)
( 1.00000000e+01,  3.75426635e-01)
( 2.00000000e+01,  7.50959854e-01)};\label{line:err_aorta_sbi_noise2}

\nextgroupplot[yticklabels={,,}, xtick={0, 5, 10, 20}, ytick={0,3,6,9,12}, ymax=7, title={$\mathrm{Re}=1000$}, ymin=0, xticklabels={,,}]
\addplot [black, solid, mark options={solid, thin}, mark=*, mark size=1.5, forget plot]
coordinates {
( 0.00000000e+00,  3.10178650e-06)
( 5.00000000e+00,  1.53091772e+00)
( 1.00000000e+01,  3.05017086e+00)
( 2.00000000e+01,  6.05798417e+00)};

\addplot [blue, solid, mark options={solid, thin}, mark=square*, mark size=1.5, forget plot]
coordinates {
( 0.00000000e+00,  1.27440281e-06)
( 5.00000000e+00,  4.59039191e-01)
( 1.00000000e+01,  9.17523856e-01)
( 2.00000000e+01,  1.83344958e+00)};

\addplot [red, solid, mark options={solid, thin}, mark=triangle*, mark size=1.5, forget plot]
coordinates {
( 0.00000000e+00,  3.07814548e-07)
( 5.00000000e+00,  2.14396861e-01)
( 1.00000000e+01,  4.27846301e-01)
( 2.00000000e+01,  8.52155239e-01)};

\nextgroupplot[xtick={0, 5, 10, 20}, ymin=40, xlabel={$\kappa$ (\%)}, ymax=100, ylabel={$e_{\mathrm{MRI}}$ (\%)}]
\addplot [black, solid, mark options={solid, thin}, mark=*, mark size=1.5, forget plot]
coordinates {
( 0.00000000e+00,  4.92727793e+01)
( 5.00000000e+00,  4.86778157e+01)
( 1.00000000e+01,  4.88207029e+01)
( 2.00000000e+01,  5.48890428e+01)};

\addplot [blue, solid, mark options={solid, thin}, mark=square*, mark size=1.5, forget plot]
coordinates {
( 0.00000000e+00,  5.00731302e+01)
( 5.00000000e+00,  5.04804451e+01)
( 1.00000000e+01,  5.13778601e+01)
( 2.00000000e+01,  5.97287299e+01)};

\addplot [red, solid, mark options={solid, thin}, mark=triangle*, mark size=1.5, forget plot]
coordinates {
( 0.00000000e+00,  5.01915016e+01)
( 5.00000000e+00,  4.53927526e+01)
( 1.00000000e+01,  4.95657434e+01)
( 2.00000000e+01,  7.96994298e+01)};

\nextgroupplot[yticklabels={,,}, xtick={0, 5, 10, 20}, ymin=40, xlabel={$\kappa$ (\%)}, ymax=100]
\addplot [black, solid, mark options={solid, thin}, mark=*, mark size=1.5, forget plot]
coordinates {
( 0.00000000e+00,  5.78198694e+01)
( 5.00000000e+00,  5.54113014e+01)
( 1.00000000e+01,  5.46254465e+01)
( 2.00000000e+01,  6.08726002e+01)};

\addplot [blue, solid, mark options={solid, thin}, mark=square*, mark size=1.5, forget plot]
coordinates {
( 0.00000000e+00,  5.75788626e+01)
( 5.00000000e+00,  5.80649285e+01)
( 1.00000000e+01,  6.15050846e+01)
( 2.00000000e+01,  7.62583799e+01)};

\addplot [red, solid, mark options={solid, thin}, mark=triangle*, mark size=1.5, forget plot]
coordinates {
( 0.00000000e+00,  5.40835180e+01)
( 5.00000000e+00,  5.41498969e+01)
( 1.00000000e+01,  6.33429779e+01)
( 2.00000000e+01,  9.99386324e+01)};

\end{groupplot}\end{tikzpicture}
 \caption{WSS reconstruction error using SBI (\textit{top row}) and MRI
 (\textit{bottom row}) as a function of noise level for various
 Reynolds numbers (\textit{columns}) and MRI grid resolutions (\textit{lines}).
 Legend: $N=3$ VPD (\ref{line:err_aorta_sbi_noise0}),
         $N=5$ VPD (\ref{line:err_aorta_sbi_noise1}),
         $N=10$ VPD (\ref{line:err_aorta_sbi_noise2}).}
 \label{fig:aorta_err_wss_noise}
\end{figure}
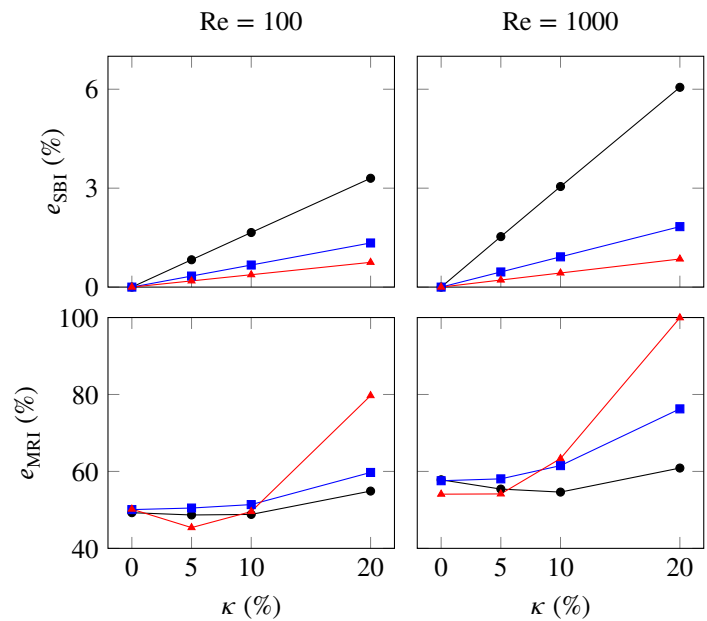

Similarly, the relationship between the WSS reconstruction error and noise
for various Reynolds numbers and MRI voxel grids
(Figure~\ref{fig:aorta_err_wss_noise}) is consistent
with findings for the stenotic vessel. For MRI, the WSS
reconstruction error increases rapidly when the noise exceeds $\kappa=5\%$
and increases most rapidly for the finest MRI voxel grid. For SBI, the WSS
reconstruction error increases linearly with the noise level with a slope
that decreases as the MRI voxel grid is refined. 

Overall, the results from the aorta test case agree with those from the
vessel, particularly the trends with respect to variations in the MRI grid
resolution and noise.

\section{Discussion}

\subsection{Summary}
This paper details a simulation-based imaging framework for blood flow
imaging and WSS reconstruction based on using numerical optimization to
fit a CFD simulation to MRI flow velocity data. The primary contribution
of the paper are two synthetic test cases---flow through a stenotic vessel
and coarctated aorta---that directly compare WSS reconstruction accuracy
using SBI and standard MRI postprocessing techniques. We found the SBI method
can accurately reconstruct WSS and is relatively insensitive to the
resolution of the MRI data and Reynolds number of the flow.
The WSS reconstruction error of SBI increases only linearly with noise in the
MRI data with a slope that is inversely proportional to the resolution of the
MRI data. Furthermore,
coarser CFD grids can be used for the SBI reconstruction than used for the
reference flow at the cost of increased sensitivity to the resolution of
the MRI data and Reynolds number of the flow. On the other hand, WSS
reconstruction from MRI data is less accurate than SBI and more sensitive to
noise, resolution of the MRI data, and the Reynolds number of the flow.

\subsection{Conclusions}
The wall shear stress reconstruction directly from MRI data
showed sensitivity to the Reynolds number of the flow, the resolution of the
MRI grid, and the noise in the MRI data. Without noise, the reconstructed
WSS converged to the true WSS distribution as the MRI grid is refined;
however, for higher noise levels the accuracy degrades as the MRI grid
is refined. This can be attributed to the noise in the domain varying more
rapidly, i.e., over a smaller length scale given by the MRI voxel spacing.
Because the MRI WSS reconstruction fits a quadratic function to this noisy
data, it is inherently sensitive to the length scale over which the noise
varies. We also observed the accuracy of the WSS reconstruction
directly from the MRI data degrades as the Reynolds number of the flow
increases for most voxel grids and noise levels considered. Lastly,
we found the error of the WSS reconstruction increases as the noise
in the MRI data increases and grows faster as the MRI grid is refined.
These findings are consistent with studies in the literature
\cite{van2013wall,potters2014measuring,potters2015volumetric,quant_wss_numerical,ha2016hemodynamic},
which show that the accuracy of WSS computations directly from MRI data
are limited by a fundamental trade-off between noise and resolution
and become less reliable as the Reynolds number of the flow increases.

On the other hand, WSS reconstruction from SBI is relatively
insensitive to the Reynolds number of the flow, the resolution
of the MRI grid, and the noise in the MRI data. For all
Reynolds numbers and noise levels considered, the WSS
reconstruction accuracy improves as the MRI grid is refined.
Furthermore, the reconstruction is reliable for all MRI grids
considered, i.e., the largest WSS reconstruction error for
the coarctated aorta test case was $e_{\mathrm{SBI}} = 6\%$, which occurred
at Reynolds number $\mathrm{Re}=1000$ with only $N=3$ VPD
and $\kappa=20\%$ noise level. The accurate WSS reconstruction
in the small data and high noise regime is attributed to the
significant amount of \textit{a priori} information leveraged
by SBI including the geometry of the domain and governing
equations (with unknown inflow boundary conditions) that
are not exploited when reconstructing WSS from MRI data alone.
Next, we observed that WSS reconstruction from SBI showed relatively
little sensitivity to the Reynolds number of the flow over the
limited range considered in this work. We expect these results
to generalize throughout the laminar regime, but break down as
the transitional and turbulent regimes are approached. Lastly,
for all Reynolds numbers and MRI grids considered, the error in
the WSS reconstruction using SBI increased linearly with noise;
the maximum error observed is $10\%$ for the stenotic vessel and
$6\%$ for the coarctated aorta which occur at a larger noise than
usually observed in practice ($\kappa=20\%$). This relatively
minor sensitivity to noise is attributed to the fact that
SBI does not directly compute WSS from a noisy field; rather,
the noisy MRI data is used to reconstruct a noise-free CFD
velocity field, which is then used to compute WSS. Because the
velocity field reconstruction is reliable in the presence of
zero-mean noise, the overall WSS reconstruction is reliable.

The SBI framework was shown to be moderately sensitive to the
resolution of the mesh used for the SBI reconstruction.
That is, a coarser mesh can be used for the SBI reconstruction
than used to simulate the true flow; however, the sensitivity
of the WSS reconstruction with respect to MRI resolution and
Reynolds number increases. This provides an opportunity to reduce
the computational cost of SBI for \textit{in vivo} applications
as it suggests there is some flexibility in designing the mesh
used for SBI provided care is taken to ensure underresolution
of the CFD velocity field is compensated with
additional MRI resolution. Multi-fidelity optimization approaches
that progressively refine inexpensive models, e.g., simulations on
coarser meshes \cite{ziems2011adaptive} or reduced-order models
\cite{zahr2015progressive}, to accelerate convergence could be
used to further reduce the cost of SBI. Alternatively, this observation
implies that a limited amount of MRI data can be used (e.g., from fast
patient scans) provided the CFD mesh used for SBI reconstruction is
relatively fine due to the insensitivity of SBI to the quality/resolution
of the MRI data in this scenario. This leads to an inherent trade-off
between scan and reconstruction time when using SBI, which is
fundamentally different and preferred than the trade-off between
scan time and reconstruction quality attributed to traditional MRI
postprocessing techniques.

\subsection{Limitations and future work}
While this study resulted in a number observations regarding
the strengths and weaknesses of using SBI to reconstruct WSS,
there are a number of limitations that offer
promising paths for future research. First, a study of the
impact of varying the point-spread function used in the
SBI reconstruction from that used to compute the true flow
would provide insight into
the sensitivity of SBI for \textit{in vivo} applications
because the point-spread function of MRI scanners is
not known with certainty, but can be modeled and estimated depending
on the MRI protocol.
Also, it would be interesting to include additional parameters
to optimize in the SBI reconstruction, e.g., material properties of
blood (with Newtonian or non-Newtonian fluid models)
and outflow conditions, to further understand the
ability of SBI to reconstruct a fully patient-specific
flow. Finally, it is important to understand how
the results obtained in this study generalize
to scenarios closer to clinical applications, e.g.,
three-dimensional domains, unsteady flows, and higher
Reynolds number flows.

\appendix

\section{Simplification of tangential component of wall traction}
\label{app:mriwss}
To derive the simplified expression for the tangential component
of the wall traction in (\ref{eqn:wss_tract_simple}), we consider
a point $x\in\partial\Omega_\mathrm{w}$.
All spatially varying quantities will be evaluated at this point;
however, for brevity, the explicit dependence on $x$ will be dropped.
Let $\Tcal$ denote the $(d-1)$-dimensional tangent space of the wall
$\partial\Omega_\mathrm{w}$ at $x$, i.e., $\Tcal$ is a linear
space such that for any $\xi\in\Tcal$, $\xi \cdot n = 0$. In addition,
let the columns of the matrix $B\in\Rbb^{d\times(d-1)}$ be an orthogonal basis
of $\Tcal$, which implies
\begin{equation} \label{eqn:defB}
 B^T n = 0, \qquad B^T B = I_{d-1}.
\end{equation}
From the no-slip condition ($v=0$ on $\partial\Omega_\mathrm{w}$),
we have that for any $\xi\in\Tcal$, $\nabla v \cdot \xi = 0$, which
implies
\begin{equation} \label{eqn:noslipB}
 \nabla v \cdot B = 0.
\end{equation}
Next, we observe that the columns of $\begin{bmatrix} B & n \end{bmatrix}$
form a basis of $\Rbb^d$ and expand the traction vector at $x$ in this
basis
\begin{equation} \label{eqn:tract_expand}
 t = B t_s + n t_n,
\end{equation}
where $t_s \in \Rbb^{d-1}$ and $t_n \in \Rbb$ are the coefficients
of the traction vector expansion. From this expansion, the tangential
component of the traction vector reduces to
\begin{equation} \label{eqn:tract_wall}
 \tau = (I-nn^T) t = B t_s,
\end{equation}
where we used (\ref{eqn:defB}) and unity of the normal vector $n$.
Furthermore, by multiplying (\ref{eqn:tract_expand}) by $B^T$ and
using (\ref{eqn:defB}), we have
\begin{equation} \label{eqn:tract_wall2}
 t_s = B^T t.
\end{equation}
Next, we substitute the expression for the traction in (\ref{eqn:wss})
into the above equation to yield
\begin{equation} \label{eqn:tract_tang}
 t_s = B^T \left[\mu(\nabla v + \nabla v^T)n + Pn\right] =
 \mu B^T (\nabla v \cdot n),
\end{equation}
where we used (\ref{eqn:defB}) and (\ref{eqn:noslipB}). Finally, we
combine (\ref{eqn:tract_wall}) and (\ref{eqn:tract_tang}) to yield
the simplified expression for the tangential component of the wall
traction in (\ref{eqn:wss_tract_simple}).

\section*{Acknowledgments}
The authors would like to thank Fritiof Hegardt for providing the aorta
mesh and simulation. This material is based upon work supported by the
Air Force Office of Scientific Research (AFOSR) under award number
FA9550-20-1-0236 and the Swedish Research Council grant 2018-03721,
the Crafoord Foundation, the Swedish strategic e-science research program
eSSENCE, and the Swedish Society of Medicine.
The content of this publication does not necessarily reflect the position
or policy of any of these supporters, and no official endorsement should be
inferred.

\bibliography{SBI}


%
%

\end{document}